\documentclass[english,fleqn,10pt]{article}
\usepackage[T1]{fontenc}
\usepackage[utf8]{inputenc}
\usepackage{morefloats}
\usepackage[letterpaper]{geometry}
\usepackage{float}
\usepackage{amsmath}
\usepackage[title]{appendix}

\usepackage{graphicx}
\usepackage{color}
\usepackage{adjustbox}
\usepackage{diagbox}
\usepackage{multirow}
\usepackage{array}
\usepackage{mathrsfs}
\usepackage{amsmath, amssymb, amsthm}
\usepackage{subfiles}
\usepackage{subcaption}
\usepackage{lipsum}
\usepackage{tikz}
\usetikzlibrary{arrows.meta,positioning,calc}

\tikzset{
  line/.style={-{Latex[length=2mm]}, thick},
  block/.style={
    draw, rounded corners, thick, align=center,
    text width=55mm, minimum height=10mm, inner sep=2mm
  },
  sideblock/.style={
    draw, rounded corners, thick, align=center,
    text width=42mm, minimum height=10mm, inner sep=2mm
  },
  sblock/.style={
    draw, rounded corners, thick, align=center,
    text width=28mm, minimum height=10mm, inner sep=2mm
  },
  note/.style={
    draw, rounded corners, thick, align=left, font=\small,
    text width=90mm, inner sep=2mm
  }
}
\usepackage{bm}
\usepackage{comment}

\usepackage[hidelinks]{hyperref}
\usepackage{url}
\usepackage[lined,boxed,linesnumbered,ruled]{algorithm2e}
\usepackage{natbib}
\usepackage{enumitem}
\setlist{nolistsep}

\newtheorem{definition}{Definition}

\usepackage{amsfonts}

\usepackage{color}
\usepackage{adjustbox}
\usepackage{diagbox}
\definecolor{black}{rgb}{0,0,0}

\definecolor{red}{rgb}{1,0,0}

\definecolor{blue}{rgb}{0,0,1}

\newcommand{\cB}{\mathcal{B}}

\newcommand{\cG}{\mathcal{G}}

\newcommand{\cL}{\mathcal{L}}

\newcommand{\cN}{\mathcal{N}}

\newcommand{\bc}{\mathbf{c}}

\newcommand{\bu}{\mathbf{u}}
\newcommand{\bv}{\mathbf{v}}

\newcommand{\bw}{\mathbf{w}}
\newcommand{\bbR}{\mathbb{R}}

\newcommand{\spec}{\text{spec}}
\newcommand{\auto}{\text{auto}}

\usepackage{multirow}

\makeatother

\usepackage{babel}
\usepackage{authblk}

\usepackage{graphicx}
\def\al#1\eal{\begin{align}#1\end{align}}
\def\als#1\eals{\begin{align*}#1\end{align*}}
\def\eq#1\eeq{\begin{eqnarray}#1\end{eqnarray}}
\def\eqs#1\eeqs{\begin{eqnarray*}#1\end{eqnarray*}}
\def\eqn#1\eeqn{\begin{equation}#1\end{equation}}
\def\eqns#1\eeqns{\begin{equation*}#1\end{equation*}}
\def\ad#1\ead{\begin{aligned}#1\end{aligned}}
\def\itm#1\eitm{\begin{itemize}#1\end{itemize}}
\def\df#1\edf{\begin{definition}#1\end{definition}}
\def\pm#1\epm{\begin{pmatrix}#1\end{pmatrix}}
\def\sl#1\esl{\begin{solution}#1\end{solution}}
\def\ex#1\eex{\begin{example}#1\end{example}}
\def\tab#1\etab{\begin{table}[!htbp] \centering#1\end{table}}
\graphicspath{{./figures/}}
\usepackage{epstopdf}
\title{
Mapped Multi-Patch Spectral Extreme Learning Machine for Partial Differential Equations
}
\author{ Yiran Wang, Suchuan Dong
}

\date{(\today)}

\begin{document}
	
\maketitle
\begin{abstract}
We study a fixed-feature solver, referred to as Spectral-ELM, in which Chebyshev--Gauss--Lobatto (CGL) differentiation matrices are applied to the nodal values of an extreme learning machine (ELM) trial function. Curved domains are treated using a mapped multi-patch formulation. Each curved quadrilateral patch is represented as the image of a reference square, physical derivatives are computed from the associated metric terms, and adjacent patches are coupled by enforcing continuity of the solution and its normal flux. Local QR orthogonalization is used to reduce near-linear dependence among the discrete features. The numerical study includes comparisons with direct CGL collocation, an analytically differentiated ELM, a global Random Feature Method, and a TransNet sampling strategy. It also includes a matched KdV test, a small three-dimensional example, and a Poisson problem on a curved domain containing a hole. For a smooth elliptic problem on a square, direct CGL collocation gives the smallest error. When applied to the same fixed-feature trial space, spectral and analytic differentiation produce comparable errors. For the curved-domain problem, mapped direct CGL collocation and mapped Spectral-ELM achieve similar accuracy with comparable numbers of unknowns. The present formulation is particularly well suited to smooth PDEs in one to three dimensions, where tensor-product grids enable accurate high-order discretizations, while sparse-grid and dimension-adaptive strategies offer promising directions for extension to higher-dimensional problems.
\end{abstract}
\par\medskip
\section{Introduction}
The numerical solution of partial differential equations depends on two closely related choices: how the unknown solution is represented and how its derivatives are evaluated. For smooth problems in one, two, or three dimensions, spectral collocation methods provide an effective answer to both questions. They represent the solution by its values at carefully selected nodes and approximate derivatives using differentiation matrices. In particular, Chebyshev--Gauss--Lobatto (CGL) collocation can attain high accuracy with relatively few grid points when the solution is sufficiently smooth.

The accuracy of spectral collocation comes from its structured polynomial approximation. This structure is also one of its restrictions. The unknowns are normally tied to all nodes of a tensor-product grid, and the geometric setting is most natural for intervals, rectangles, and rectangular boxes. Although coordinate mappings and domain decomposition can extend spectral methods to more general domains, the approximation remains closely connected to the underlying nodal grid. These features motivate the search for alternative representations that retain accurate numerical differentiation while allowing the solution to be described by a smaller set of coefficients.

Neural-network PDE solvers approach the representation problem from a different direction. Rather than assigning an independent unknown to each node, they approximate the solution by a parameterized function. Physics-informed neural networks determine the network parameters by minimizing residuals associated with the governing equation, boundary and initial conditions, and available data \cite{raissi2017physics,karniadakis2021physics}. This functional representation is flexible, but the resulting optimization problem is nonlinear and can be sensitive to the network architecture, residual sampling, loss weights, initialization, and training algorithm \cite{cuomo2022scientific,cyr2020robust}. The cost and variability of training can therefore become as important as the approximation properties of the network itself.

Fixed-feature methods retain the functional representation of a neural network while avoiding the optimization of all network parameters. In an Extreme Learning Machine (ELM), the hidden-layer weights and biases are selected in advance and remain fixed, while only the output coefficients are determined from the governing equation \cite{huang2004extreme,huang2006extreme}. For linear PDEs, this usually leads to a linear least-squares problem. For nonlinear PDEs, the unknowns are still limited to the output coefficients, so the resulting nonlinear system is substantially smaller than a fully trained network problem. The main task is no longer to train a large collection of parameters, but to construct a useful feature space and determine the coefficients within that space.

This viewpoint has led to several related approaches. Local ELM methods use domain decomposition and impose coupling conditions between neighboring subdomains \cite{dong2021local}, while constrained trial functions incorporate boundary information directly into the approximation \cite{schiassi2021extreme}. The Random Feature Method combines local randomized features with collocation, equation scaling, and domain decomposition \cite{chen2022bridging}. TransNet distributes neuron hyperplanes throughout the domain and selects a common scale so that the resulting feature space can be reused across related approximation problems \cite{zhang2023transnet}. The Subspace Neural Network first constructs a collection of neural basis functions and then solves the PDE within their span \cite{xu2025subspace}. Despite their differences, these methods all separate the construction of a finite-dimensional trial space from the final calculation of the PDE solution.

Once a fixed-feature trial space has been constructed, its derivatives still need to be evaluated. The usual choice is to differentiate the features analytically or through automatic differentiation. This is natural when the feature functions and differential operators are simple. However, it ties the implementation to the selected activation function and requires the derivative calculation to be developed for each operator and derivative order. The effort becomes more noticeable for equations involving several high-order derivatives or for mapped formulations in which reference-coordinate derivatives must be combined with geometric metric terms.

Classical spectral differentiation offers another possibility. Instead of differentiating each feature separately, one may first evaluate the fixed-feature approximation at CGL nodes and then apply the standard CGL differentiation matrices to these nodal values. This leads to the Spectral-ELM studied in this paper. The solution is represented by the output coefficients of a fixed ELM, but its derivatives are calculated through the same nodal operators used in Chebyshev spectral collocation.

This construction brings together two useful properties. The fixed-feature representation allows the number of unknown coefficients to be chosen independently of the number of points at which the PDE residual is evaluated. A relatively small set of coefficients can therefore define the trial solution on a denser CGL grid. At the same time, spectral differentiation provides a common procedure for computing derivatives of different orders after the trial function has been sampled. The differentiation stage does not need to be redesigned when the feature family changes, provided that the feature values can be evaluated on the grid.

The combination is not automatically advantageous, and understanding its behavior requires separating two different sources of error. The first is the restriction imposed by the fixed-feature trial space. Direct CGL collocation allows the solution value at every node to vary independently, whereas Spectral-ELM restricts the nodal solution to the space generated by the selected features. This restriction reduces the number of coefficient unknowns, but it may also exclude nodal solutions that are available to the direct spectral method. The approximation quality therefore depends on how well the sampled features represent the solution.

The second source of error is the differentiation procedure. Analytic differentiation acts directly on the feature functions, while CGL differentiation acts on the polynomial interpolant determined by their nodal values. The two approaches use the same trial space but need not produce identical derivatives. Spectral differentiation is attractive because it provides one reusable nodal mechanism, but its accuracy depends on whether the CGL grid resolves the variation of the features. Analytic differentiation avoids this interpolation step, although it requires feature-specific derivative formulas.

These observations lead naturally to two comparisons. Direct CGL collocation provides the appropriate reference for assessing the effect of replacing a full nodal representation with a fixed-feature representation. A matched analytic-derivative ELM provides the appropriate reference for assessing the effect of the derivative calculation itself. In the latter comparison, the trial function, collocation points, feature realization, coefficient unknowns, residual equations, and algebraic solver must remain unchanged. Only the evaluation of derivatives should differ. This separation allows the roles of representation and differentiation to be examined independently.

The same distinction is important for nonlinear and time-dependent PDEs. Small changes in the trial function or in the treatment of initial and boundary conditions can affect the nonlinear algebraic problem and make derivative comparisons difficult to interpret. For this reason, the Korteweg--de Vries experiment in this paper uses an identical fixed-feature trial function for spectral and analytic differentiation. The two versions also share the same collocation points, random features, unknown coefficients, residual formulation, initialization, and nonlinear solver. The comparison is designed to reveal the effect of the derivative calculation rather than differences in the underlying approximation.

A further question concerns geometry. CGL differentiation matrices are most naturally constructed on tensor-product reference domains. Simply introducing them into an ELM formulation would therefore restrict the method to regular geometries and would give up much of the geometric flexibility associated with feature-based collocation. To make spectral differentiation useful in a broader setting, a mechanism is needed to connect the structured reference grid to curved physical domains.

We address this issue through a mapped multi-patch construction. A curved domain is divided into quadrilateral patches, each obtained by mapping a reference square into the physical domain. The CGL points, fixed features, and differentiation matrices are defined on the reference patch. Physical derivatives are then calculated through the metric information of the map. In this way, the same reference domain formulation can be used on patches with different sizes, orientations, and curved boundaries.

The patch maps are constructed from their boundary curves using transfinite interpolation. The underlying geometric ideas are classical and are closely related to curvilinear coordinates and mapped spectral elements \cite{patera1984spectral,gordon1973curvilinear}. Their role here is to provide a direct connection between the reference CGL grid and the physical domain. Adjacent patches are coupled by requiring continuity of the solution and continuity of the physical normal flux across their common interfaces. The same construction can be arranged around an interior boundary, allowing the method to treat a curved domain containing a hole.

The fixed-feature representation also introduces a conditioning issue. Randomly sampled features may be distinct as continuous functions but nearly dependent when evaluated on a finite nodal grid. Such dependence can make the coefficient calculation unstable and can obscure the actual approximation ability of the feature space. We therefore construct a local QR basis from the sampled features on each patch. This step produces a numerically stable representation of the resolved feature space before the PDE, boundary, and interface conditions are assembled into the least-squares system.

These ingredients define the main purpose of the paper. We develop a fixed-feature PDE solver in which CGL matrices provide the derivative approximation, local QR bases stabilize the discrete feature spaces, and mapped patches extend the construction to curved domains. The method is studied not as a replacement for classical spectral collocation, but as a way to combine a reduced feature representation with established high-order differentiation tools. This viewpoint makes it possible to examine when the fixed-feature restriction is effective and when the full nodal spectral space remains preferable.

The numerical study is organized around this question. A smooth elliptic problem on a square is used to compare Spectral-ELM with direct CGL collocation and with an ELM based on analytic feature derivatives. A global Random Feature Method and a TransNet sampling strategy provide additional fixed-feature controls. These experiments separate the influence of the feature space from that of the differentiation procedure and also show how the results vary across different random feature realizations.

The remaining examples examine other parts of the formulation. The matched Korteweg--de Vries problem tests the two derivative strategies for a nonlinear equation containing a third-order derivative. A small three-dimensional problem confirms that the tensor-product construction extends beyond two dimensions. Finally, a Poisson equation on a curved multi-patch domain with an interior hole tests the geometric mapping, transformed differential operators, curved boundary conditions, and interface coupling.

The experiments show that direct CGL collocation remains the most accurate method for the smooth elliptic problem on a square. This result reflects the strength of the full nodal polynomial space on a problem for which classical spectral collocation is already well suited. The fixed-feature formulations use fewer coefficient unknowns, but their accuracy depends on the quality and dimension of the sampled feature space. The local QR construction improves the numerical representation of this space, although it does not remove the approximation error caused by an insufficient set of features.

When the trial space is held fixed, spectral and analytic differentiation give comparable errors in the tests considered here. Analytic differentiation is slightly more accurate in the matched Korteweg--de Vries experiment at the tested resolution, while spectral differentiation provides a single nodal procedure for evaluating the required derivatives. These results suggest that the choice of feature space often has a greater effect on the final PDE error than the choice between the two differentiation strategies.

The curved-domain experiment illustrates the role of the proposed formulation more clearly. For the Poisson problem with an interior hole, mapped Spectral-ELM and mapped direct CGL collocation achieve similar accuracy when compared at comparable numbers of unknowns. The benefit of Spectral-ELM in this setting is not an improvement over spectral collocation on a simple square. It is the ability to combine a compact local feature space, standard reference-domain differentiation, and physically consistent patch coupling within one formulation.

The present method is intended for smooth PDEs in one to three dimensions, where tensor-product CGL grids remain computationally practical. Although the feature representation can reduce the number of unknown coefficients, the differentiation matrices still require the trial function to be evaluated on the tensor grid. The method is therefore not designed as a high-dimensional PDE solver. Sparse grids, low-rank tensor approximations, dimension-adaptive constructions, and non-tensor nodal sets offer possible routes for extending the same fixed-feature and numerical-differentiation idea to higher-dimensional problems.

The remainder of the paper is organized as follows. Section 2 presents the Spectral-ELM formulation, the local QR construction, and the mapped multi-patch treatment of curved domains. Section 3 gives the controlled comparisons and the equation-specific numerical experiments. Section 4 summarizes the main findings, discusses the scope of the method, and considers possible extensions.

\section{Spectral Extreme Learning Machine}

Consider the following problem:
\begin{subequations}
    \begin{align}
    \cL u+\cN u &= f, \quad x\in \Omega,\\ 
     \cB u &= g,\quad x\in \partial \Omega,
    \end{align}
\end{subequations}
where $\cL$ and $\cN$ are linear and nonlinear operators, respectively; $\cB$ is the boundary operator; $\Omega\subset\bbR^d$, with $d\leq 3$; and $f$ and $g$ are the source term and boundary data.

\subsection{Spectral Collocation Discretization}
A spectral method represents a function by an expansion in a set of basis functions. In a collocation formulation, derivatives of the nodal interpolant are evaluated by differentiation matrices and the PDE is enforced at the nodes. Smooth functions can converge rapidly as the number of nodes increases. This accuracy applies to the resolved nodal interpolant; using an ELM trial space adds a separate feature-approximation error and therefore does not guarantee higher accuracy.

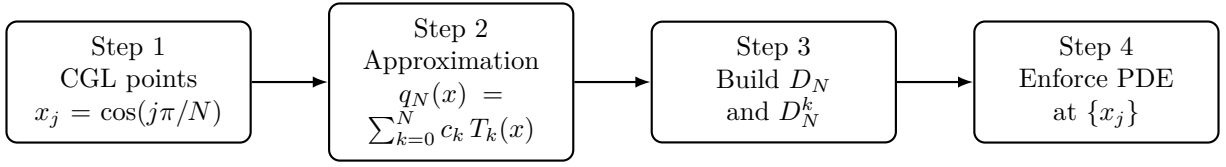
\begin{figure}[!htbp]
\centering
\resizebox{0.98\linewidth}{!}{%
\begin{tikzpicture}[node distance=10mm and 10mm]

\node[sblock] (s1) {Step 1\\CGL points\\ $x_j=\cos(j\pi/N)$};
\node[sblock, right=of s1] (s2) {Step 2\\Approximation\\ $q_N(x)=\sum_{k=0}^N c_k\,T_k(x)$};
\node[sblock, right=of s2] (s3) {Step 3\\Build $D_N$\\and $D_N^k$};
\node[sblock, right=of s3] (s4) {Step 4\\Enforce PDE\\at $\{x_j\}$};

\draw[line] (s1) -- (s2);
\draw[line] (s2) -- (s3);
\draw[line] (s3) -- (s4);

\end{tikzpicture}
}
\caption{Spectral collocation workflow used in this paper: CGL points,
Chebyshev interpolation, differentiation matrix construction, and pointwise
enforcement of the governing equation. For the illustrative one-dimensional
Poisson problem, the last step assembles
$-D_N^2\mathbf{u}=\mathbf{f}$ and imposes $u_0=u_N=0$; higher derivatives are
evaluated as $D_N^k=(D_N)^k$, with
$(D_N^k\mathbf{v})_j=q^{(k)}(x_j)$.}
\label{fig:spectral_collocation_steps}
\end{figure}

The spectral collocation method approximates the PDE solution by enforcing the equation at a set of collocation points. The process has four steps: generate Gauss--Lobatto points, approximate the PDE solution with Chebyshev polynomials, construct the collocation differentiation matrices, and enforce the PDE at the collocation points.
To clarify the implementation pipeline, Figure~\ref{fig:spectral_collocation_steps} summarizes the four main steps of the spectral collocation method used in this work, from constructing Chebyshev--Gauss--Lobatto points to enforcing the governing equations at collocation nodes.
This workflow will be used repeatedly in the subsequent sections to construct spectral differentiation matrices and enforce PDE residuals at collocation points within the proposed Spectral-ELM framework.

\textbf{Chebyshev--Gauss--Lobatto Points.} When Chebyshev--Gauss--Lobatto (CGL) points are employed, the method benefits from the superior approximation properties of Chebyshev polynomials and the inclusion of boundary points, which facilitates the treatment of boundary conditions. This makes the method particularly well suited to problems with smooth solutions.

The Chebyshev--Gauss--Lobatto points are defined by
\als
x_j = \cos\left(\frac{j\pi}{N}\right), \quad j=0,\ldots, N,
\eals
where $N$ is the number of subintervals (or the degree of the approximating polynomial). These points are the extrema of the Chebyshev polynomial of the first kind, $T_N(x)$, and include the endpoints of the interval. Their clustering near the endpoints helps resolve boundary layers or rapid variations near the edges of the domain. The inclusion of boundary points allows boundary conditions to be imposed directly, while the nonuniform spacing mitigates interpolation errors such as Runge's phenomenon.

\textbf{Polynomial Approximation of the Solution.} The solution $u(x)$ to the PDE is approximated by a polynomial interpolant $q_N(x)$, which passes through the values of $u(x)$ at the CGL points. The interpolant is expressed as:
\als
q_N(x) = \sum\limits_{k=0}^{N} c_k T_k(x),
\eals
where $T_k(x)$ are the Chebyshev polynomials of the first kind, and $c_k$ are the coefficients determined by the function values $u(x_j)$ at the CGL points. The accuracy of the interpolation improves significantly with an increasing number of points $N$, provided the underlying solution is smooth. 

\textbf{Construction of the Differentiation Matrix.} To compute the derivative of the solution at the CGL points, a differentiation matrix $D_N$ is constructed. Let $\bv=[v_0, v_1,\ldots, v_N]^T$ represent the values of the function $v(x)$ at the CGL points. The derivative $\bw = [w_0,w_1,\ldots,w_N]^T$ at these points is approximated as
 \begin{align*}
     D_N\left[\begin{array}{c}v_0 \\ \vdots \\ v_N\end{array}\right]=\left[\begin{array}{c}w_0 \\ \vdots \\ w_N\end{array}\right],
 \end{align*}
 where $v_0, \ldots, v_N$ are the values of $v$ and $w_0,\ldots,w_N$
 approximate its derivatives at the CGL points. The implicit definition of $D_N$ is as follows:
 \begin{definition}
     Suppose $q_N$ is the degree-$N$ polynomial interpolant of $v$ at the CGL points. The associated Chebyshev differentiation matrix $D_N$ satisfies the following properties:
     \begin{enumerate}
	        \item Interpolate $v$ by $q_N$ such that $q_N(x_j)=v(x_j)$;
 	        \item Differentiate the interpolant at the grid points $x_j$:
 	        \begin{equation}
		            w_j = (D_N\bv)_j = q_N'(x_j). \label{first_order_mat}
 		        \end{equation}
 	    \end{enumerate}
 \end{definition}

 Based on this definition, $D_N$ has the following entries when $x_0, \ldots, x_N$ form a CGL grid:
 \begin{align*}
     &(D_N)_{00} = \frac{2N^2+1}{6}, \quad (D_N)_{NN}= -\frac{2N^2+1}{6}, \\
     &(D_N)_{jj} = \frac{-x_j}{2(1-x_j^2)}, \text{ for } 1\leq j\leq N-1,\\
     &(D_N)_{ij} =  \frac{c_i}{c_j}\frac{(-1)^{i+j}}{x_i-x_j}, \text{ for }i\neq j,
 \end{align*}
 where $c_0=c_N=2$ and $c_j=1$ for $1\leq j\leq N-1$. For higher-order derivatives, we define $D_N^k=(D_N)^k$ as the matrix associated with the $k$th derivative. In particular, corresponding to \eqref{first_order_mat},
 \begin{equation}
     (D_N^k\bv)_j = q_N^{(k)}(x_j). \label{k_order_mat}
 \end{equation}
\textbf{Enforcing the Governing Equations.} The original PDE is enforced at the CGL points. For instance, consider the one-dimensional Poisson equation
\al
-u_{xx}(x)&=f(x),\quad x\in [-1,1],\\
u(-1)&=u(1)=0,
\label{1d_poisson}
\eal
the process is as follows. First, the differentiation matrix $D_N^2$ is constructed for the second derivative. Next, the function $f(x)$ is evaluated at the CGL points, yielding $\mathbf{f}$. The resulting system,  
\als
-D_N^2 \bu = \mathbf{f},
\eals
is solved for $\bu$, subject to the boundary conditions $u_0=u_N=0$. The discrete solution $\bu$ provides approximate values of $u(x)$ at the CGL points, and the polynomial interpolant $q_N(x)$ can be used to recover the solution over the entire domain.

This formulation provides a robust and efficient framework for solving PDEs with smooth solutions while ensuring high accuracy and straightforward implementation of boundary conditions.

\subsection{Extreme Learning Machine Method}
An Extreme Learning Machine (ELM) uses a randomized feedforward neural network and has gained significant attention for its computational efficiency and simplicity. It is characterized by randomly initialized hidden-layer parameters and a closed-form solution for the output-layer weights. This architecture makes an ELM particularly suitable for solving PDEs because it can approximate solutions while avoiding the computational cost of iterative optimization over all network parameters.

An ELM represents the solution of a PDE as an approximate function of the form:
\als
u(x) \approx \mathcal{G}(x;\theta) =\sum\limits_{i=1}^L\beta_i \sigma(w_i\cdot  x+b_i),
\eals
where $\sigma$ is the activation function, $w_i$ and $b_i$ are the weight and bias of the $i$th hidden neuron, $\beta_i$ are the trainable output weights, $L$ is the number of hidden neurons, and $x\in \Omega$ represents input variables such as spatial and temporal coordinates. The hidden-layer parameters $w_i$ and $b_i$ are randomly initialized and fixed, while only the output weights $\beta_i$ are optimized during training.

The method enforces the governing equation at selected collocation points, transforming the continuous differential equation into an algebraic system. This system is then solved to determine the output weights $\beta_i$.

\textbf{Methodology.} The ELM solution process involves the following steps. First, the domain is discretized into collocation points $\{x_j\}_{j=1}^{N_c}$ and, if applicable, boundary points $\{x_k\}_{k=1}^{N_b}$. At these points, the governing PDE and boundary conditions are enforced. The ELM approximation $\cG(x;\theta)$ is substituted into the PDE, converting the differential equation $\cL[u(x)]=f(x)$ into a set of algebraic equations. For each collocation point $x_j$, the following equation is imposed:
\als
\cL[\cG(x_j;\theta)] =f(x_j).
\eals
This step yields $N_c$ algebraic equations. Similarly, the boundary conditions $\cB[\cG(x_k;\theta)]=g(x_k)$ are enforced at the boundary points $\{x_k\}$, adding $N_b$ equations.

The enforcement of these equations results in the following algebraic system:
\als
\mathbf{A} \beta = \mathbf{g},
\eals
where $\mathbf{A}$ contains the activation functions and their derivatives evaluated at the collocation points, together with coefficients from the differential operator $\cL$; $\beta = [\beta_1,\beta_2,\ldots, \beta_L]^T$ is the vector of unknown output weights; and $\mathbf{g}$ contains the source and boundary data.

\begin{figure}[!htbp]
\centering
\begin{tikzpicture}[node distance=8mm]

\tikzstyle{every node}=[font=\footnotesize]

\node[block] (ans)
{ELM ansatz\\
$\mathcal{G}(x;\theta)=\sum_{i=1}^M \beta_i\sigma(w_i\cdot x+b_i)$\\
$\{w_i,b_i\}$ fixed; $\beta$ unknown};

\node[block, below=of ans] (pts)
{Select discrete points\\
collocation $\{x_j\}_{j=1}^{N_c}$\\
boundary $\{x_k\}_{k=1}^{N_b}$};

\node[block, below=of pts] (op)
{Enforce PDE and BC\\
$\mathcal{L}[\mathcal{G}(x_j;\theta)]=f(x_j)$\\
$\mathcal{B}[\mathcal{G}(x_k;\theta)]=g(x_k)$};

\node[block, below=of op] (sys)
{Assemble algebraic system\\
$\mathbf{A}\beta=\mathbf{g}$};

\node[block, below=of sys] (solve)
{Solve for $\beta$\\
(linear system / least squares)};

\draw[line] (ans) -- (pts);
\draw[line] (pts) -- (op);
\draw[line] (op) -- (sys);
\draw[line] (sys) -- (solve);

\end{tikzpicture}
\caption{ELM-based PDE solution mechanism: discretize the domain, enforce the governing equations at selected points, and solve a linear system for the output weights.}
\label{fig:elm_pde_flow}
\end{figure}
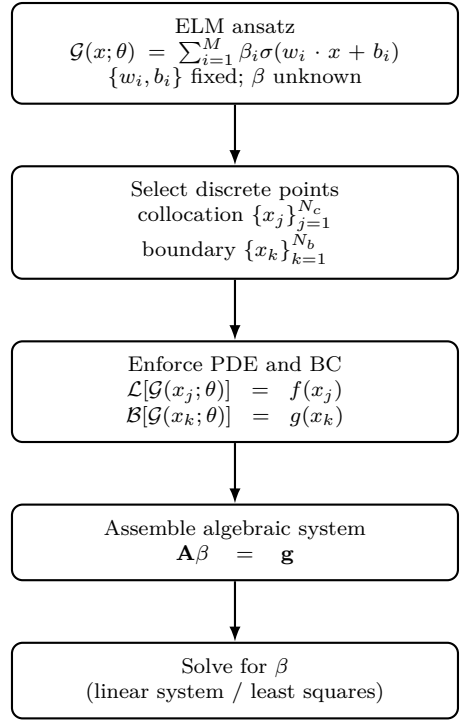

The output weights $\beta_i$ are determined by solving this computationally efficient linear system. As illustrated in Figure~\ref{fig:elm_pde_flow}, the ELM-based PDE solver transforms the continuous differential equation into a finite-dimensional linear algebraic system by enforcing the governing equations and boundary conditions at discrete collocation points. The hidden-layer parameters are randomly initialized and fixed, so the only unknowns are the output weights $\beta$. Consequently, training reduces to solving a structured linear system rather than performing iterative nonlinear optimization. The computational complexity is therefore primarily determined by the construction of the system matrix $\mathbf{A}$ and the solution of the resulting least-squares problem.

In practical implementations, the entries of $\mathbf{A}$ depend on the evaluation of the activation functions and their derivatives at the collocation points. In standard ELM-based PDE solvers, these derivatives are typically obtained analytically or via automatic differentiation of the neural representation. This derivative evaluation step plays a central role in the overall accuracy and stability of the method, and will be revisited in the subsequent subsection when introducing the Spectral-ELM formulation.

\textbf{Example: Solving the Poisson Equation.} Consider the one-dimensional Poisson equation on the interval $[0,1]$:
\al
-u_{xx}(x) = f(x), \quad x\in [0,1], \quad u(0)=u(1)=0. \label{poisson_1d_elm}
\eal
The domain $[0,1]$ is discretized into $N_c$ collocation points $\{x_j\}$ and the boundary points $x=0$ and $x=1$. Substituting the ELM approximation $u(x)\approx \sum_{i=1}^{L} \beta_i \sigma(w_i\cdot x+b_i)$ into the equation transforms the Poisson equation into a set of algebraic equations at the collocation points. The second derivative $-u_{xx}$ is computed analytically from the ELM representation by differentiating the activation functions. The boundary conditions $u(0)=0$ and $u(1)=0$ are imposed directly in the algebraic system. Solving the resulting equations yields the weights $\beta_i$, which define the approximate solution.

We next give a matrix-vector representation for computing solution derivatives, using \eqref{poisson_1d_elm} as an example. The goal is to represent the second derivative $u_{xx}$ by differentiating the hidden-layer outputs. These derivatives are efficiently calculated through automatic differentiation, which directly uses the computational graph of the ELM.

First, an ELM approximates the solution $u(x)$ by the neural-network output $U(x)$:
\als
U(x) = \sum\limits_{i=1}^{M} \beta_i\sigma(w_ix +b_i),
\eals
where $\sigma(\cdot)$ is the activation function, $w_i,b_i\in \bbR$ are the weights and biases of the $i$-th neuron, and $\beta_i\in \bbR$ are the trainable output weights. The hidden layer outputs are represented as 
\als
\mathbf{h}(x)=\left[\begin{array}{c}
	\sigma\left(w_1 \cdot x+b_1\right) \\
	\sigma\left(w_2 \cdot x+b_2\right) \\
	\vdots \\
	\sigma\left(w_M \cdot x+b_M\right)
\end{array}\right],
\eals 
and the solution can be written compactly as $U(x)=\mathbf{\beta}^T \mathbf{h}(x)$, where $\mathbf{\beta} = [\beta_1,\beta_2,\ldots, \beta_M]^T$. To compute the second derivative $U_{xx}(x)$, the second derivatives of the activation functions, weighted by the squared neuron weights, must be determined and stored. Automatic differentiation calculates these derivatives directly from the computational graph, enabling a straightforward and efficient implementation.

Next, the second-order derivative $U_{xx}$ is expressed as 
\als
U_{xx}(x) = \sum\limits_{i=1}^{M} \beta_i \sigma''(w_ix+b_i)w_i^2,
\eals
where $\sigma''(\cdot)$ represents the second derivative of the activation function and $w_i^2$ is the squared weight of the $i$th neuron. Evaluating this derivative at the collocation points $x_1,\ldots, x_N$, we construct a matrix $D^{(2)}$ whose $i$th column contains the second derivatives of the $i$th neuron's output. The entries of $D^{(2)}$ are defined as
\als
(D^{(2)})_{ji}=\sigma''(w_ix_j+b_i)w_i^2,
\eals
where $x_j$ represents the $j$-th collocation point, and $i$ indexes the neurons in the hidden layer. The matrix $D^{(2)}$ has dimensions $N\times M$, where $N$ is the number of collocation points and $M$ is the number of hidden-layer neurons.

Finally, the second derivative of $U(x)$ at all collocation points is computed as
\al
\mathbf{U}_{xx}=D^{(2)}\beta, \label{2matrix-vec}
\eal
where $\mathbf{U}_{xx}=[U_{xx}(x_1),U_{xx}(x_2),\ldots, U_{xx}(x_N)]$. The matrix $D^{(2)}$ is constructed by applying automatic differentiation to the activation functions in the computational graph, thereby evaluating their second derivatives at the collocation points.

Using automatic differentiation to compute the derivatives ensures accuracy and eliminates the need to derive and implement analytical derivatives for each activation function. This approach is advantageous when the activation function is complex or when extending the model to higher-order derivatives. For example, consider the common activation function $\sigma(z) = \tanh(z)$. Its first and second derivatives are $\sigma'(z)=1-\tanh^2(z)$ and $\sigma''(z)=-2\tanh(z)(1-\tanh^2(z))$, respectively. By leveraging automatic differentiation, these derivatives are computed automatically.

Once the second derivative $U_{xx}$ is computed in \eqref{2matrix-vec}, we can transform the continuous differential equation \eqref{poisson_1d_elm} into a system of algebraic equations:
\als
-D^{(2)}\beta=\mathbf{f},
\eals
where $\mathbf{f}$ represents the values of $f(x)$ at the collocation points. This formulation enables computation of the second derivatives while maintaining the simplicity and scalability of the ELM framework.
\subsection{Spectral Extreme Learning Machine Method}
\begin{figure}[!htbp]
\centering
\begin{tikzpicture}[node distance=10mm and 14mm]

\node[note] (key)
{\textbf{Key idea.} Compute derivatives using spectral differentiation matrices on CGL points (instead of automatic differentiation), and solve only for the output weights $\beta$ via a linear/least-squares system.};

\node[block, below=of key] (pde)
{$\mathcal{L}u+\mathcal{N}u=f$ in $\Omega$, \quad $\mathcal{B}u=g$ on $\partial\Omega$};

\node[block, below=of pde] (grid)
{Choose CGL points \ $\{x_j\}_{j=0}^{N}$};

\node[block, below=of grid] (elm)
{ELM ansatz \ $u(x)\approx \mathcal{G}(x;\theta)=\sum_{i=1}^M \beta_i\sigma(w_i\cdot x+b_i)$\\
$\theta=\{w_i,b_i\}$ fixed; solve for $\beta$};

\node[block, below=of elm] (eval)
{Evaluate on grid: \ $\mathbf{U}=[\mathcal{G}(x_j;\theta)]_{j=0}^{N}$};

\node[block, below=of eval] (enf)
{Enforce PDE/BC at collocation points \ $\Rightarrow$ residual equations};

\node[block, below=of enf] (ls)
{Assemble linear system / LS solve: \ $\mathbf{A}\beta=\mathbf{g}$};

\node[block, below=of ls] (out)
{Obtain $\beta$ and reconstruct solution $u(x)$ (interpolant / evaluation)};

\node[sideblock, right=of eval, xshift=8mm] (spec)
{Build spectral matrices\\ $D_N,\ D_N^2,\ \ldots$};

\node[sideblock, below=of spec] (deriv)
{Compute derivatives on grid\\ $\mathbf{U}_x=D_N\mathbf{U}$,\\
$\mathbf{U}_{xx}=D_N^2\mathbf{U}$, etc.};

\draw[line] (key) -- (pde);
\draw[line] (pde) -- (grid);
\draw[line] (grid) -- (elm);
\draw[line] (elm) -- (eval);
\draw[line] (eval) -- (enf);
\draw[line] (enf) -- (ls);
\draw[line] (ls) -- (out);

\draw[line] (grid.east) |- (spec.west);
\draw[line] (eval.east) -- (spec.west);
\draw[line] (spec) -- (deriv);
\draw[line] (deriv.west) |- (enf.east);

\end{tikzpicture}
\caption{Overall workflow of the Spectral-ELM method.}
\label{fig:overall_spectral_elm}
\end{figure}
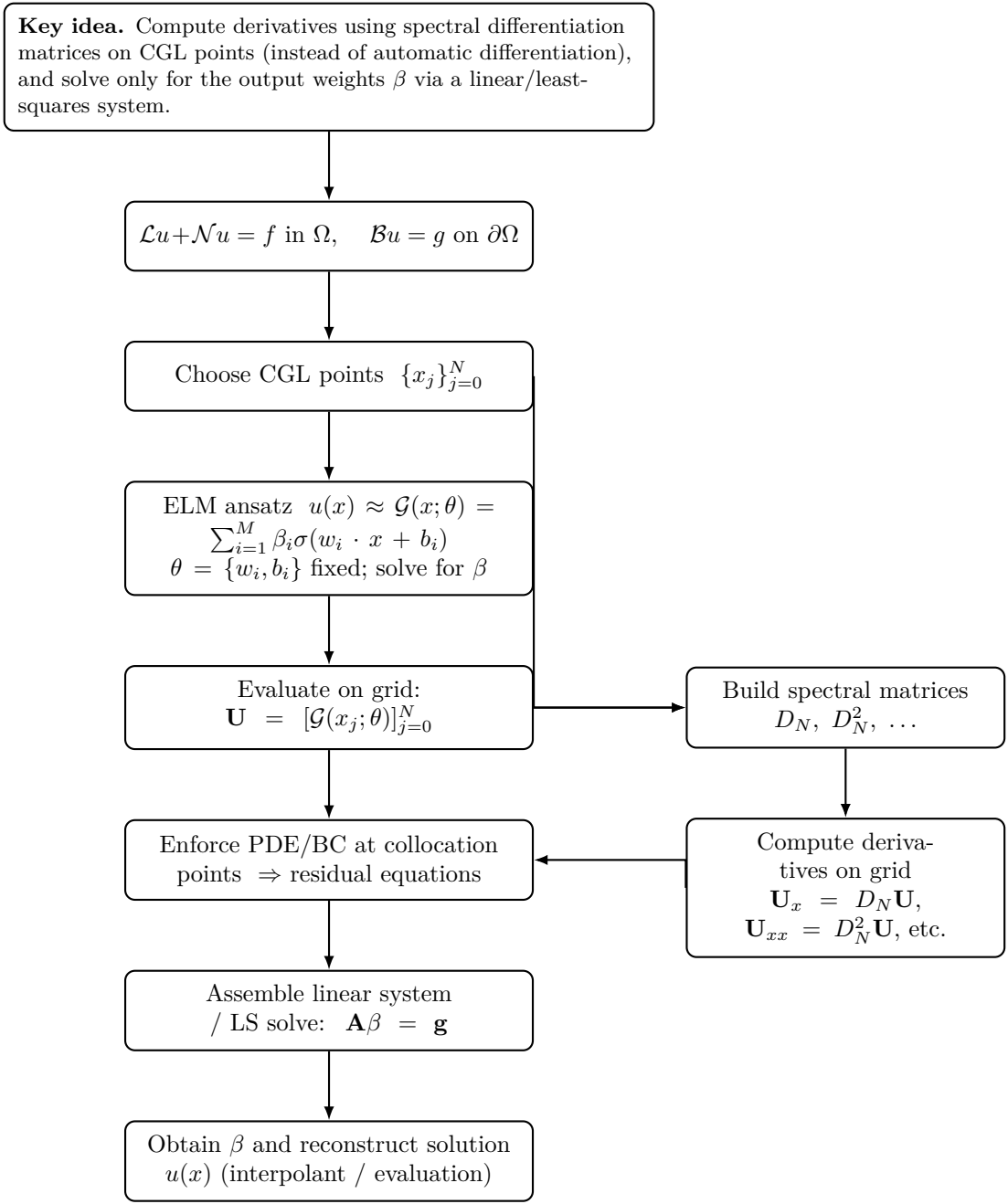

The Spectral-ELM method combines the approximation power of an ELM with the high accuracy of spectral discretization. As illustrated in Figure~\ref{fig:overall_spectral_elm}, the method begins by discretizing the domain using Chebyshev--Gauss--Lobatto (CGL) points, at which the governing equations and boundary conditions are enforced. The solution is represented by a fixed random-feature ELM ansatz in which only the output weights are unknown. Instead of differentiating the neural representation through automatic differentiation, spectral differentiation matrices constructed at the CGL points are used to compute spatial and, when applicable, temporal derivatives. These derivatives are then used to assemble the residual equations at the collocation points, leading to a structured linear system for the output weights. Unlike standard ELM implementations, Spectral-ELM decouples function approximation from derivative computation, allowing derivative accuracy to be governed by spectral theory rather than neural smoothness. This combination of spectral collocation and linear ELM training yields a computationally efficient and numerically precise framework for solving linear and nonlinear PDEs.

Consider the following general PDE problem:
\als
\cL u +\cN u&= f, \quad x \in \Omega,\\
\cB u &= g,\quad x\in \partial \Omega,
\eals
where $\cL$ and $\cN$ are linear and nonlinear operators, respectively, $\cB$ represents the boundary operator, and $f$ and $g$ are the source term and boundary data. The domain $\Omega\subset \bbR^d$ is assumed to be bounded, with $d\leq 3$. The Spectral-ELM framework approximates $u(x)$ by integrating an ELM representation with the spectral collocation method.

\textbf{Methodology.} In the Spectral-ELM framework, the solution $u(x)$ is approximated as
\als
u(x)\approx \cG(x;\theta) =\sum\limits_{i=1}^M \beta_i\sigma(w_ix+b_i),
\eals
where $\sigma(\cdot)$ is the activation function, $w_i\in \bbR^d$ and $b_i\in \bbR$ are the fixed weight and bias of the $i$th neuron, $\beta_i$ are the trainable output weights, and $M$ is the number of hidden neurons. Unlike traditional neural networks, Spectral-ELM randomly initializes and fixes $w_i$ and $b_i$, while optimizing only the output weights $\beta_i$ to satisfy the PDE.

\textbf{Spectral Collocation for Derivative Computation.} Spectral-ELM uses CGL points to compute derivatives by spectral collocation. These points are defined as
\als
x_j = \cos\left(\frac{j\pi}{N}\right), \quad j=0,1,\ldots, N,
\eals
where $N$ is the number of intervals in the discretized domain. The clustering of points near the boundaries facilitates accurate resolution of steep gradients and simplifies the enforcement of boundary conditions.

The derivatives of the approximation $U(x)$ are computed using spectral differentiation matrices. Let $\mathbf{U}=[U(x_0),U(x_1),\ldots,U(x_N)]^T$ represent the solution values at the CGL points. The first and second derivatives of $U(x)$ are approximated as
\als
 \mathbf{U}_x=D_N\mathbf{U},\quad \mathbf{U}_{xx}=D_N^2\mathbf{U},
\eals
where $D_N$ is the first-order spectral differentiation matrix, and $D_N^2=D_N\cdot D_N$ is the second-order differentiation matrix.

To solve the PDE, the neural representation $\cG(x;\theta)$ is substituted into the governing equations. The spectral differentiation matrices $D_N$ and $D_N^2$ are used to compute the derivatives of $\cG(x;\theta)$ at the CGL points. For example, for the one-dimensional Poisson equation \eqref{poisson_1d_elm}, the discrete form is written as:
\als
-D_N^2 \mathbf{U}=\mathbf{f},
\eals
where $\mathbf{U}$ is the vector of approximate solution values at the CGL points and $\mathbf{f}$ is the source term evaluated at these points. The boundary conditions are enforced directly by modifying the system matrix to impose $u(-1)=u(1)=0$. The output weights $\beta=[\beta_1,\beta_2,\ldots,\beta_M]^T$ are determined by solving a linear system that enforces the PDE and boundary conditions at the collocation points:
\als
\mathbf{A}\beta =\mathbf{g},
\eals
where $\mathbf{A}$ is constructed using the neural representation of $\cG(x;\theta)$ and the spectral differentiation matrices, and $\mathbf{g}$ contains the source terms and boundary values.

The Spectral-ELM method applies a standard nodal differentiation operator to a compact fixed-feature trial space. It avoids differentiating the activation functions, but it is accurate only when the CGL grid resolves their nodal variation. Whether this is better than analytic differentiation or direct nodal collocation must be checked for each problem.

\subsection{Mapped Multi-Patch Extension for Complex Geometries}
\label{sec:mapped_multipatch}

The global CGL construction above is limited to an interval or a
tensor-product domain. To handle a complex geometry, we divide the physical
domain into curved quadrilateral patches,
\begin{equation}
	\overline{\Omega}=\bigcup_{e=1}^{N_e}\overline{\Omega}_e,
	\qquad
	\Omega_e\cap\Omega_{e'}=\varnothing \quad(e\ne e').
\end{equation}
Each physical patch is mapped from the reference square:
\begin{equation}
	F_e:[-1,1]^2\longrightarrow\Omega_e,
	\qquad
	(x,y)=F_e(\xi,\eta).
	\label{eq:patch_map}
\end{equation}
The maps may be given analytically, constructed by transfinite interpolation,
or obtained from a high-order mesh. This follows the standard mapped spectral
element idea \cite{patera1984spectral,gordon1973curvilinear}. The complex
boundary is represented by the patch maps, while the same CGL grid is used on
every reference square.

\paragraph{How the form of $F_e$ is chosen.}
The map is determined from the geometry before the PDE is solved. First choose
four physical curves for patch $e$: the left and right curves
$C_L(\eta),C_R(\eta)$ and the bottom and top curves
$C_B(\xi),C_T(\xi)$. The curves must meet at the same four corners. A direct
choice is the Gordon--Hall transfinite map
\begin{align}
	F_e(\xi,\eta)={}&\frac{1-\xi}{2}C_L(\eta)
	+\frac{1+\xi}{2}C_R(\eta)
	+\frac{1-\eta}{2}C_B(\xi)
	+\frac{1+\eta}{2}C_T(\xi) \notag\\
	&-\frac{(1-\xi)(1-\eta)}{4}P_{LB}
	-\frac{(1-\xi)(1+\eta)}{4}P_{LT} \notag\\
	&-\frac{(1+\xi)(1-\eta)}{4}P_{RB}
	-\frac{(1+\xi)(1+\eta)}{4}P_{RT},
	\label{eq:gordon_hall_map}
\end{align}
where $P_{LB},P_{LT},P_{RB},P_{RT}$ are the four corner points. The first
line blends the four curves into the patch. The second and third lines remove
the corner contributions that were counted twice. Equation
\eqref{eq:gordon_hall_map} is the default choice when only the four boundary
curves are known. An analytic map can be used instead when the geometry has a
natural coordinate system.

A proposed map is accepted only if its Jacobian $J_e$ does not vanish and has
one sign throughout the reference square. This is checked on a grid that is
finer than the collocation grid. If $J_e$ is zero or changes sign, the map has
folded. The remedy is geometric: split the physical region into smaller
patches, change the interior patch curves, or use a high-order mesh map. The
PDE solver cannot repair a folded $F_e$.

For the wavy annulus in Section~\ref{sec:wavy_annulus}, the natural map is even
simpler. In angular patch $e$, let
\begin{equation}
	\theta_e(\eta)=\theta_{e,\mathrm{mid}}
	+\frac{\Delta\theta_e}{2}\eta,
	\qquad s(\xi)=\frac{1+\xi}{2},
\end{equation}
and set
\begin{equation}
	r_e(\xi,\eta)=[1-s(\xi)]r_{\rm in}(\theta_e)
	+s(\xi)r_{\rm out}(\theta_e),
	\qquad
	F_e(\xi,\eta)=
	\begin{bmatrix}r_e\cos\theta_e\\r_e\sin\theta_e\end{bmatrix}.
	\label{eq:wavy_annulus_map}
\end{equation}
Here $\xi=-1$ is exactly the inner boundary, $\xi=1$ is exactly the outer
boundary, and $\eta=-1$ and $\eta=1$ are the two radial interfaces. Adjacent angular
patches use the same curve at their shared angle, so their interface nodes
match automatically.

The reference and physical calculations should not be confused. The CGL
nodes, random features, and matrices $D_\xi,D_\eta$ are built on the regular
square. The map sends those nodes to physical coordinates. The source and
boundary data are evaluated at the physical coordinates; the map derivatives
give the metric terms; and the PDE and physical normal flux are assembled with
those metrics. After the coefficients are solved, evaluating the solution at a
new physical point means locating its patch, finding the corresponding
$(\xi,\eta)$ (analytically or by a two-variable Newton solve), and evaluating
the local feature expansion there. The differential equation is therefore
enforced in physical space, even though the reusable grid and basis data live
on the reference square.

Let $D$ be the one-dimensional CGL differentiation matrix with $n$ nodes. We
define the two reference differentiation matrices by
\begin{equation}
	D_\xi=D\otimes I_n,
	\qquad
	D_\eta=I_n\otimes D.
\end{equation}
At the tensor-product nodes, the derivatives of the map give
\begin{equation}
	J_e=x_\xi y_\eta-x_\eta y_\xi,
	\quad
	\xi_x=\frac{y_\eta}{J_e},\quad
	\xi_y=-\frac{x_\eta}{J_e},\quad
	\eta_x=-\frac{y_\xi}{J_e},\quad
	\eta_y=\frac{x_\xi}{J_e}.
	\label{eq:inverse_metrics}
\end{equation}
The physical first-derivative matrices on patch $e$ are therefore
\begin{align}
	D_x^{(e)}
	&=\operatorname{diag}(\xi_x)D_\xi
	+\operatorname{diag}(\eta_x)D_\eta,\\
	D_y^{(e)}
	&=\operatorname{diag}(\xi_y)D_\xi
	+\operatorname{diag}(\eta_y)D_\eta.
	\label{eq:mapped_first_derivatives}
\end{align}

For the Poisson operator, define
\begin{equation}
	g^{\xi\xi}=\xi_x^2+\xi_y^2,\qquad
	g^{\xi\eta}=\xi_x\eta_x+\xi_y\eta_y,\qquad
	g^{\eta\eta}=\eta_x^2+\eta_y^2.
\end{equation}
The mapped discrete Laplacian is written in conservative form as
\begin{align}
	L_e={}&\operatorname{diag}(J_e)^{-1}
	\Big[
	D_\xi\operatorname{diag}(J_e g^{\xi\xi})D_\xi
	+D_\xi\operatorname{diag}(J_e g^{\xi\eta})D_\eta \notag\\
	&\hspace{22mm}
	+D_\eta\operatorname{diag}(J_e g^{\xi\eta})D_\xi
	+D_\eta\operatorname{diag}(J_e g^{\eta\eta})D_\eta
	\Big].
	\label{eq:mapped_laplacian}
\end{align}
Thus, changing the geometry only changes the map and its metric terms; the
reference CGL differentiation matrices remain unchanged.

On patch $e$, we use a local ELM representation
\begin{equation}
	u_e\big(F_e(\xi,\eta)\big)
	=\sum_{m=1}^{M_e}\beta_{e,m}
	\sigma\!\left(w_{e,m}^{T}[\xi,\eta]^T+b_{e,m}\right).
	\label{eq:local_elm_patch}
\end{equation}
Let $H_e\in\mathbb{R}^{n^2\times M_e}$ contain the fixed features evaluated at
the local CGL nodes. Then
\begin{equation}
	\mathbf u_e=H_e\boldsymbol\beta_e,
	\qquad
	\Delta\mathbf u_e=L_eH_e\boldsymbol\beta_e.
\end{equation}
The use of local random features is related to domain-decomposed ELM and local
randomized neural methods \cite{dong2021local,sun2024local,lee2025nonoverlapping}.
The difference here is that the derivatives are computed with mapped CGL
matrices instead of differentiating the neural features.

For an interface
$\Gamma_{ee'}=\partial\Omega_e\cap\partial\Omega_{e'}$, we impose continuity
of the solution and normal flux:
\begin{align}
	H_e^{\Gamma}\boldsymbol\beta_e
	-P_{ee'}H_{e'}^{\Gamma}\boldsymbol\beta_{e'}&=0,
	\label{eq:interface_c0}\\
	N_eH_e\boldsymbol\beta_e
	+P_{ee'}N_{e'}H_{e'}\boldsymbol\beta_{e'}&=0.
	\label{eq:interface_flux}
\end{align}
Here, $P_{ee'}$ matches the nodes on the two sides of the interface and
$N_e=n_xD_x^{(e)}+n_yD_y^{(e)}$ is the outward normal derivative. Boundary
conditions are imposed at CGL points on the physical boundary. The PDE,
boundary, and interface equations are combined into one least-squares system
for all local output coefficients.

The raw random-feature matrices may be ill-conditioned. For each patch, we
therefore compute a thin QR factorization
\begin{equation}
	H_e=Q_eR_e,
	\qquad
	\boldsymbol\gamma_e=R_e\boldsymbol\beta_e,
\end{equation}
and solve with the equivalent basis $Q_e\boldsymbol\gamma_e$. This changes the
basis but not the discrete trial space when $H_e$ has full column rank. We use
the same row scaling obtained from the nodal spectral operator for both the
direct spectral and Spectral-ELM systems, followed by column equilibration.

This construction handles curved boundaries, holes, and nonconvex domains when
a valid patch decomposition is available. It is a boundary-fitted high-order
method rather than a fully meshless method. It also does not remove the $n^d$
growth of a tensor grid inside one $d$-dimensional patch. Therefore, geometric
complexity and high-dimensional scaling remain separate issues.

\section{Numerical Examples}
In this section, we test the methods on regular domains and on a curved domain with a hole. We denote the architecture of the neural network by $M_{\text{arch}}=[m_{\text{in}},M,m_{\text{out}}]$, where $m_{\text{in}}$ and $m_{\text{out}}$ are the numbers of input and output variables. In all tests, $m_{\text{out}}=1$. When the problem is time-dependent, $m_{\text{in}}=d+1$, where $d$ is the spatial dimension. Otherwise, $m_{\text{in}}=d$. The hidden weights and biases are sampled from $[-R_m,R_m]$, and the hyperbolic tangent is used in all examples. We use $N_x$ for the number of collocation points in each spatial direction and $N_t$ for the number of time points. The curved-domain example uses local mapped grids and states its parameters separately.

In the following numerical results, we use the maximum error
\als
e=\max|u_{\text{approx}}-u_{\text{ref}}|,
\eals
where $u_{\text{approx}}$ and $u_{\text{ref}}$ are the approximate and reference solutions. We further define $e_{\spec}$ and $e_{\auto}$ as the maximum errors obtained using spectral collocation and automatic differentiation, respectively.

We implement two differentiation strategies to formulate the system of equations. In the automatic differentiation (AD) approach, forward-mode AD computes exact derivatives of the neural features. The spectral collocation approach instead uses differentiation matrices, such as the Chebyshev $D_{xx}$ matrix, to approximate derivatives on a specified grid. To find the optimal coefficients $\beta$, we minimize the residual $\|r\|_2$ using the SciPy trust-region reflective (TRF) nonlinear least-squares solver. We apply a strict stopping tolerance of $10^{-12}$, with a maximum of 400 function evaluations. To ensure numerical stability, all feature blocks are column-scaled to unit root-mean-square (RMS) magnitude over the stacked PDE and boundary-condition rows. Furthermore, the residuals are weighted by quadrature weights $w$, such as those from the Clenshaw--Curtis or Gauss--Lobatto--Jacobi rules.

\subsection{Controlled elliptic baselines and dimensional scope}
\label{sec:controlled_baselines}

This test keeps the PDE and evaluation points fixed while changing the trial space or the derivative calculation. On $\Omega=(-1,1)^2$, we solve
\begin{equation}
	-\Delta u+u=f, \qquad u=0\quad\text{on }\partial\Omega,
\end{equation}
with the manufactured solution
\begin{equation}
	u(x,y)=(1-x^2)(1-y^2)\exp\{\sin(2x+0.5y)\}.
\end{equation}
All errors are evaluated at the same 20,000 off-grid points. Direct CGL uses $25^2=625$ nodal unknowns. Each feature method uses 300 features. Spectral-ELM and analytic-derivative ELM use the same Gaussian features, the same boundary factor $(1-x^2)(1-y^2)$, the same CGL points, and the same local QR basis. The global RFM control uses uniformly sampled tanh features, strong-form residual equations, boundary equations, and equation-wise scaling. The TransNet control samples unit feature directions and hyperplane offsets and uses the domain-scaled value $\gamma=1$. Five fixed random seeds are used for every randomized feature method. For each run, the maximum error is $\max_{1\leq j\leq 20{,}000}|u_h(\mathbf{x}_j)-u(\mathbf{x}_j)|$. The ``median max. error'' is the median of these five seed-specific maximum errors, and the ``error range'' is their minimum-to-maximum interval. The ``median condition'' is the median over the five seeds of the matrix 2-norm condition number after the column scaling used in the coefficient solve. Direct CGL is deterministic, so only its single maximum error and the condition number of its nodal system are reported.

\begin{table}[H]
	\centering
	\small
	\begin{tabular}{lrrrr}
		\hline
		Method & Unknowns & Median max. error & Error range & Median condition \\
		\hline
		Direct CGL & 625 & $3.09\times10^{-10}$ & deterministic & $9.66\times10^4$ \\
		Spectral-ELM & 300 & $4.03\times10^{-8}$ & $[1.77,8.83]\times10^{-8}$ & $4.44\times10^1$ \\
		Analytic ELM & 300 & $8.66\times10^{-8}$ & $[2.45,14.9]\times10^{-8}$ & $4.60\times10^1$ \\
		Global RFM & 300 & $5.75\times10^{-5}$ & $[2.24,13.4]\times10^{-5}$ & $6.74\times10^9$ \\
		TransNet sampling & 300 & $1.97\times10^{-7}$ & $[1.80,4.52]\times10^{-7}$ & $3.14\times10^1$ \\
		\hline
	\end{tabular}
	\caption{Controlled two-dimensional elliptic comparison. For each randomized
		feature method, the table reports the median of five seed-specific maximum
		errors on the common 20,000-point test set, the corresponding
		$[\min,\max]$ error interval, and the median matrix 2-norm condition number of
		the column-scaled coefficient system. Direct CGL is deterministic; its single
		maximum error and nodal-system condition number are shown, and no seed range
		applies. The RFM and TransNet rows have the limited scope stated in
		Section~1.}
	\label{tab:controlled_elliptic}
\end{table}

\begin{figure}[H]
	\centering
	\includegraphics[width=0.78\linewidth]{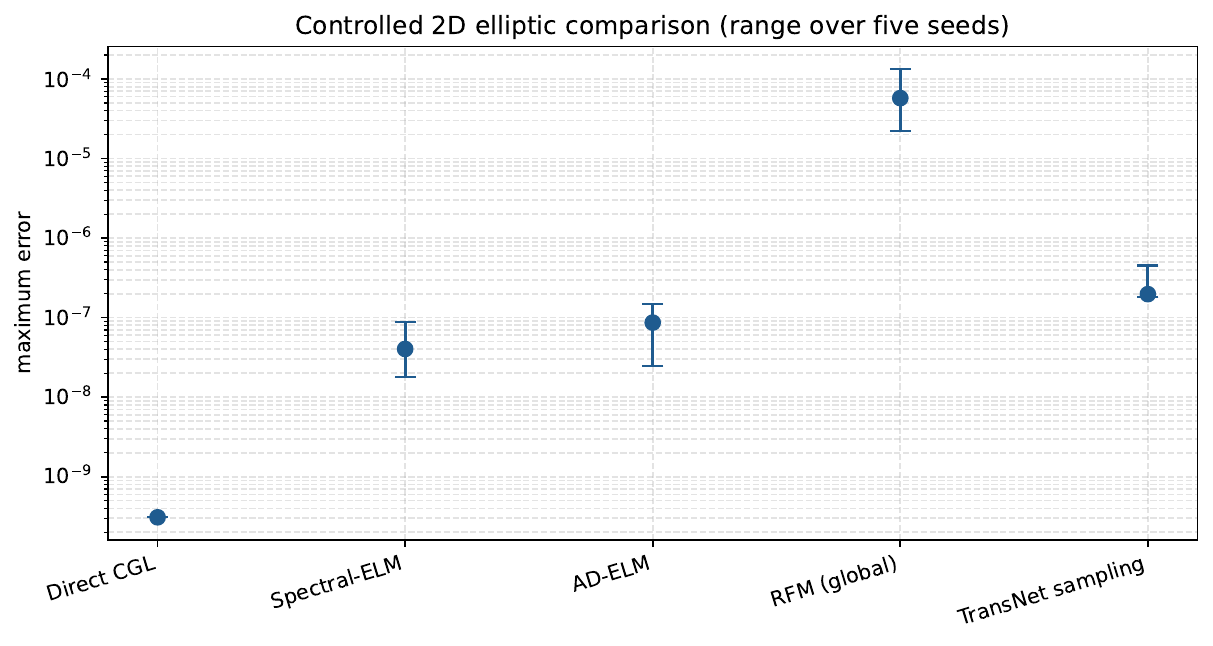}
	\caption{Maximum error for the controlled elliptic test. A point is the median; a bar gives the minimum and maximum over five seeds. Direct CGL is deterministic.}
	\label{fig:controlled_elliptic}
\end{figure}

The direct nodal method is the most accurate method in Table~\ref{tab:controlled_elliptic}. This is expected: the solution is smooth, the geometry is a square, and direct CGL uses the complete nodal polynomial space. The feature methods compress that space to 300 coefficients and introduce a feature-selection error. Spectral-ELM and analytic ELM are within the same order of accuracy, so this example does not show a universal advantage for either derivative calculation. The simplified RFM result also shows that row scaling alone does not remove near dependence in a global random feature matrix. The TransNet sampling result is much better than the global RFM result, but it is still a feature-sampling control rather than the fully tuned transferable method.

We also repeat the manufactured problem in three dimensions with
$u=\prod_{j=1}^3(1-x_j^2)\exp\{\sin(2x_1+0.5x_2-0.7x_3)\}$.
Direct CGL with $11^3=1331$ unknowns gives a maximum error of
$1.78\times10^{-3}$. With 250 features, Spectral-ELM gives a median maximum
error of $5.93\times10^{-3}$ over three seeds, with range
$[3.55\times10^{-3},1.17\times10^{-2}]$. This test confirms that the formula
extends to three coordinates, but it also displays the limitation: increasing
$n$ by one direction changes $n^2$ nodes into $n^3$ nodes. We therefore make no
claim of scalability to high dimension.

\subsection{Linear Poisson equation}
\textbf{Problem setup.} We consider the following Poisson equation in the domain $\Omega = [-1,1]$. 
\begin{eqnarray*}
    \begin{aligned}
        -u_{xx}=f, \quad x\in \Omega,\\
        u=g,\quad x\in \partial\Omega.
    \end{aligned}
    \label{eqn:linear_poisson}
\end{eqnarray*}
This system has the exact solution $u=\sin(x+0.05)$. For the test, we partition $\Omega$ into $n_e=2$ elements. Within each element, we distribute $N_x=40$ Gauss--Lobatto--Jacobi (GLJ) collocation nodes. To maintain a global solution across the subdomains, we explicitly enforce $C^0$ and $C^1$ continuity at the interface by appending interfacial constraints to the least-squares system. This ensures that both the solution $u$ and its first derivative $u_x$ remain continuous across the $n_e-1$ interfaces.

\begin{figure}[!tb]
	\centerline{
		\includegraphics[width=0.45\textwidth]{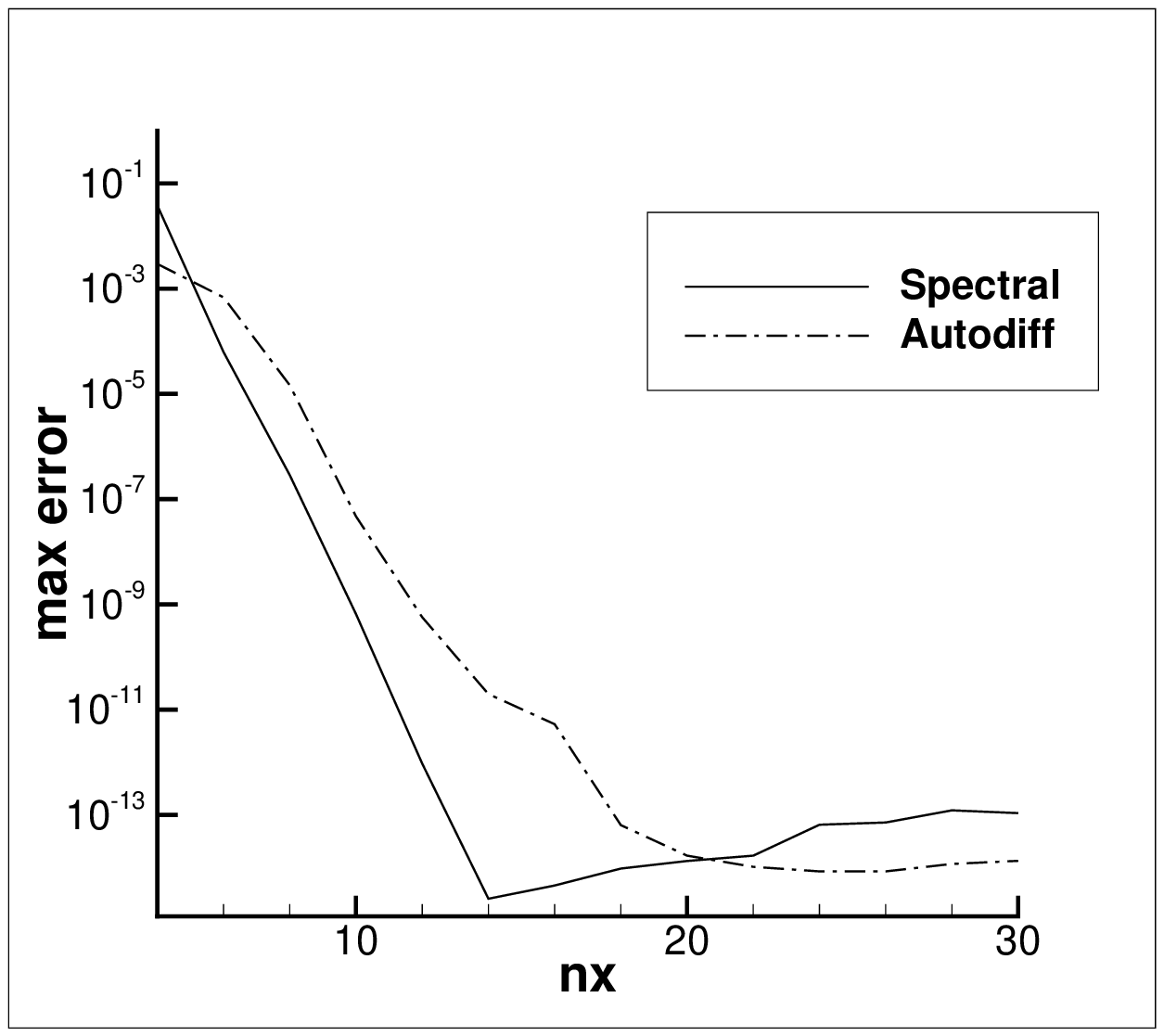}(a)
		\includegraphics[width=0.45\textwidth]{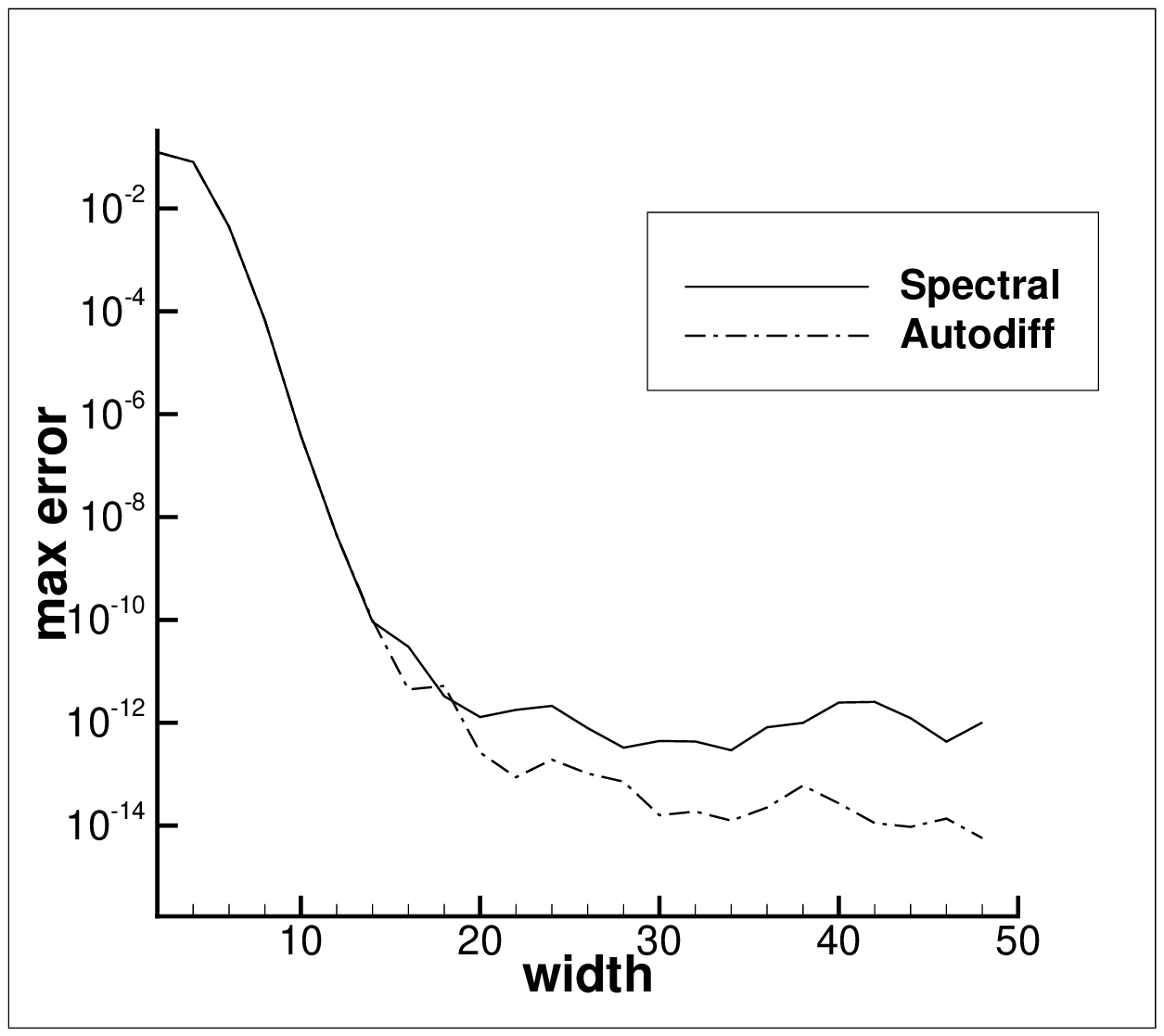}(b)
	}
		\caption{Maximum-norm approximation errors from spectral collocation and automatic differentiation versus the number of collocation points $N_x$ (left) and hidden-layer width $M$ (right). Left: $R_m=1$ and $M=40$; right: $R_m=1$ and $N_x=40$.
		}\label{fig:linear_poi_compare_max}
\end{figure}

\textbf{Numerical results.} The sensitivity of the approximation to the random magnitude $R_m$ is summarized in Table~\ref{tab:rm_linear_poi}. We observe a high-accuracy regime for $R_m$ between 0.5 and 4.0, where the AD variant achieves a maximum error $e_{\text{auto}}$ near machine precision ($9.99\times 10^{-15}$), while the spectral variant maintains a stable error $e_{\text{spec}}$ of approximately $2.07\times10^{-11}$.
\begin{table}[!htbp]
\centering
\begin{tabular}{ccccccc}
	\hline
	$R_m$ & 1E-3& 0.5 &1&2&4&8 \\ \hline 
	$e_{\text{spec}}$ &1.28E-2&2.32E-11& 2.07E-11& 3.78E-11&2.08E-11&2.18E-10\\ 
	$e_{\text{auto}}$ &1.51E-4& 1.00E-13& 3.92E-14& 9.99E-15&5.06E-14&1.86E-11\\ 
	\hline
\end{tabular}
\caption{Comparison of approximation errors from spectral collocation and automatic differentiation for different random magnitudes $R_m$. We use $M=40$ and $N_x=40$. The quantities $e_{\spec}$ and $e_{\auto}$ are the maximum errors obtained using spectral collocation and automatic differentiation, respectively.
}
\label{tab:rm_linear_poi}
\end{table}
The convergence characteristics relative to the collocation density $N_x$ and hidden-layer width $M$ are illustrated in Figure~\ref{fig:linear_poi_compare_max}. Panel (a) demonstrates exponential convergence with respect to the number of collocation points $N_x$ in the range $[14,20]$, which is characteristic of spectral-type methods. Panel (b) shows that the error decreases consistently as the hidden-layer width $M$ increases from 10 to 100 neurons, indicating that the dimensionality of the random feature space is a primary factor in the accuracy of the neural approximation.

\begin{figure}[!tb]
	\centerline{
		\includegraphics[width=0.45\textwidth]{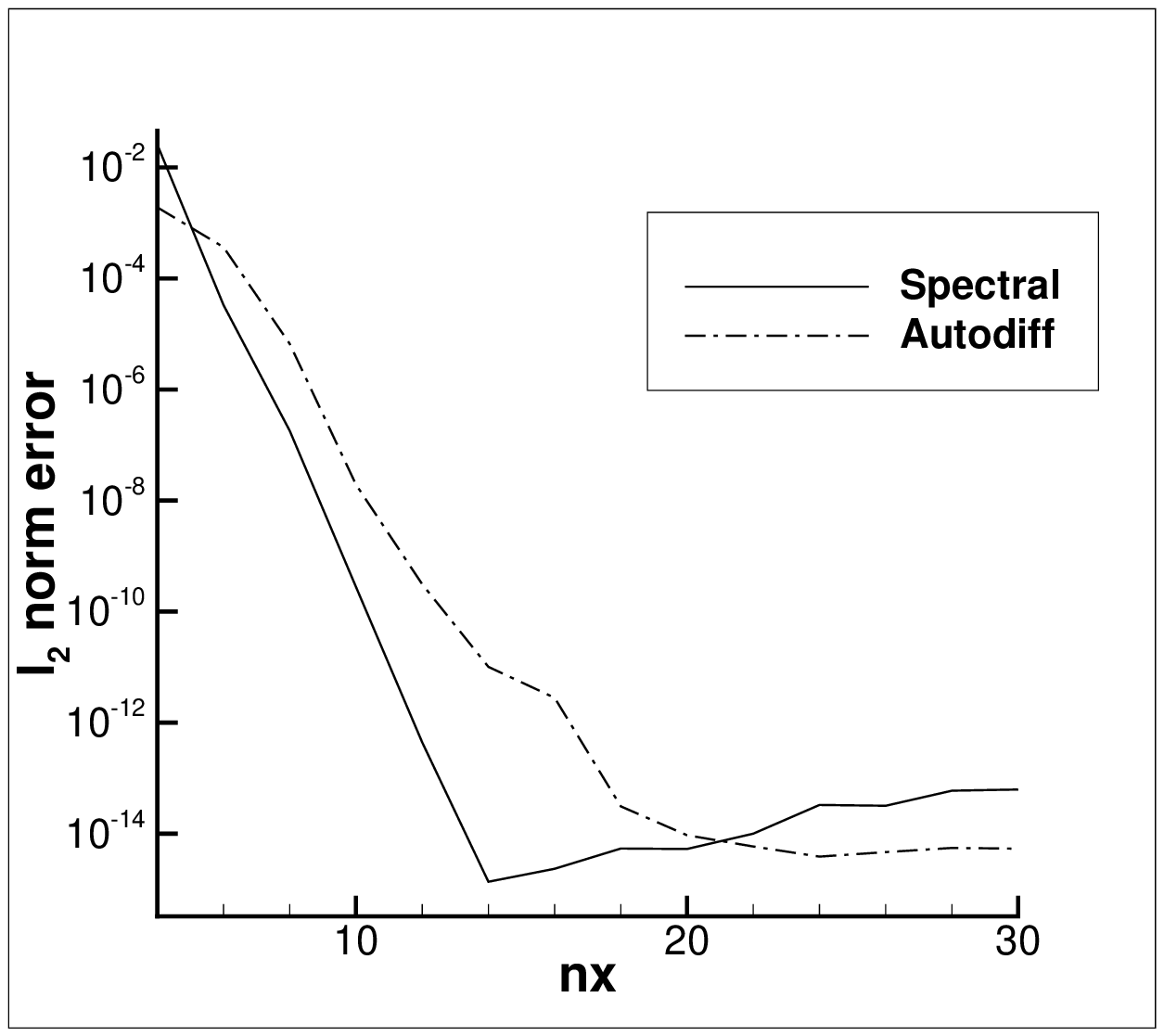}(a)
		\includegraphics[width=0.45\textwidth]{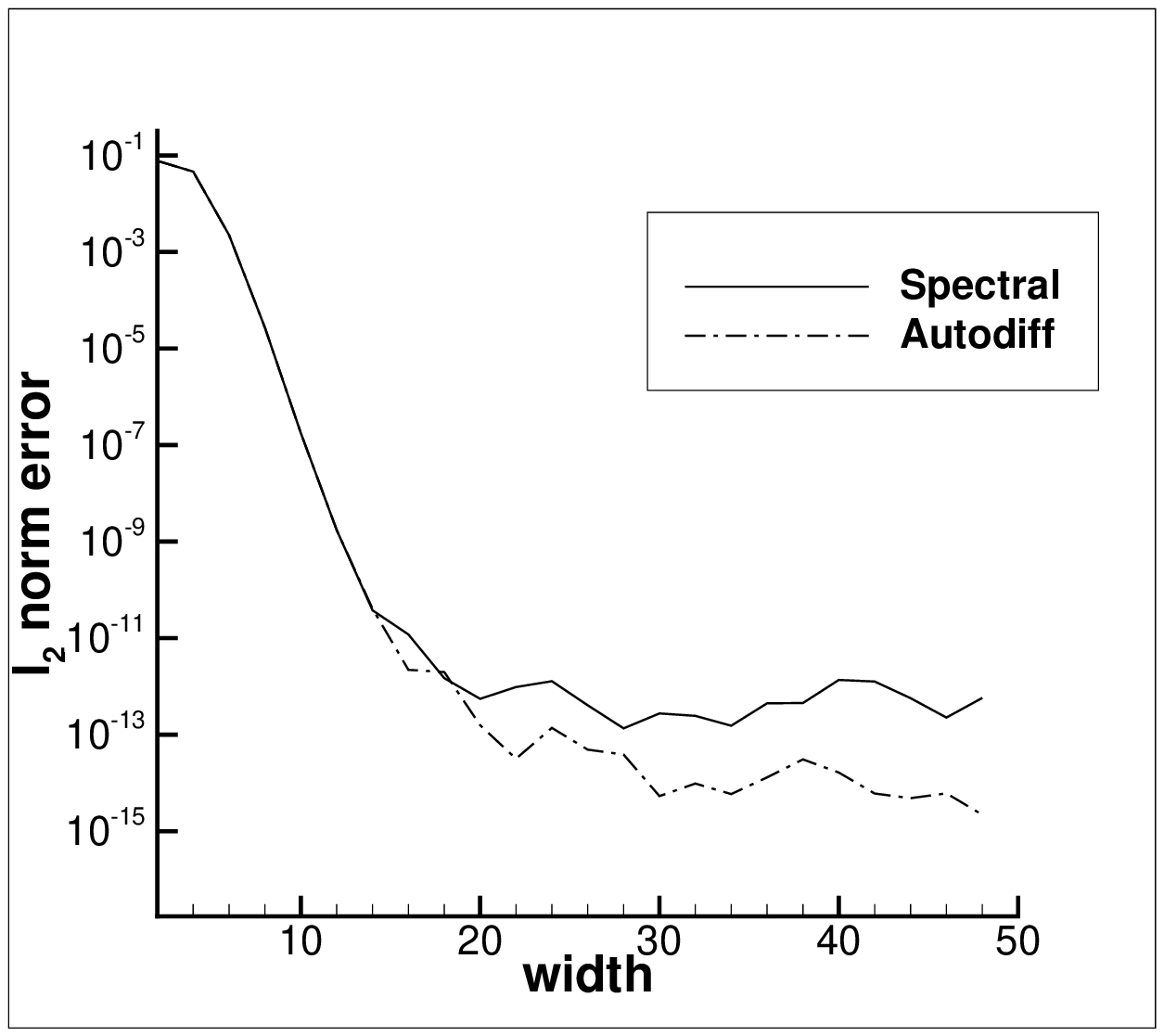}(b)
	}
		\caption{$L^2$-norm approximation errors from spectral collocation and automatic differentiation versus the number of collocation points $N_x$ (left) and hidden-layer width $M$ (right). Left: $R_m=1$ and $M=40$; right: $R_m=1$ and $N_x=40$.
		}\label{fig:linear_poi_compare_l2}
\end{figure}
\subsection{Nonlinear Poisson equation}
\textbf{Problem setup.} We consider the nonlinear Poisson equation on the two-dimensional domain $[-1,1]\times[-1,1]$:
\begin{eqnarray*}
	\begin{aligned}
		-\Delta u-u^3=f, \quad x\in \Omega,\\
		u=g,\quad x\in \partial\Omega.
	\end{aligned}
	\label{eqn:nonlinear_poisson}
\end{eqnarray*}
We use $u(x,y)=\sin(x)\sin(y)$ as the analytical solution. This case introduces significant complexity due to the cubic nonlinearity, requiring the iterative solution of the resulting discrete system.

For the two-dimensional problem, we distribute a grid of $N_x\times N_y$ GLJ collocation nodes, with $N_x=N_y$ ranging from 4 to 50. Given the increased dimensionality and nonlinearity, the hidden-layer width $M$ is typically scaled relative to the number of collocation points; in our tests, we vary $M$ from 20 to 300 neurons. The nonlinear system is solved using the TRF algorithm. To facilitate convergence and mitigate the risk of falling into local minima, we use a tight stopping tolerance of $10^{-12}$ and column scaling.

\begin{figure}[!tb]
	\centerline{
\includegraphics[width=0.45\textwidth]{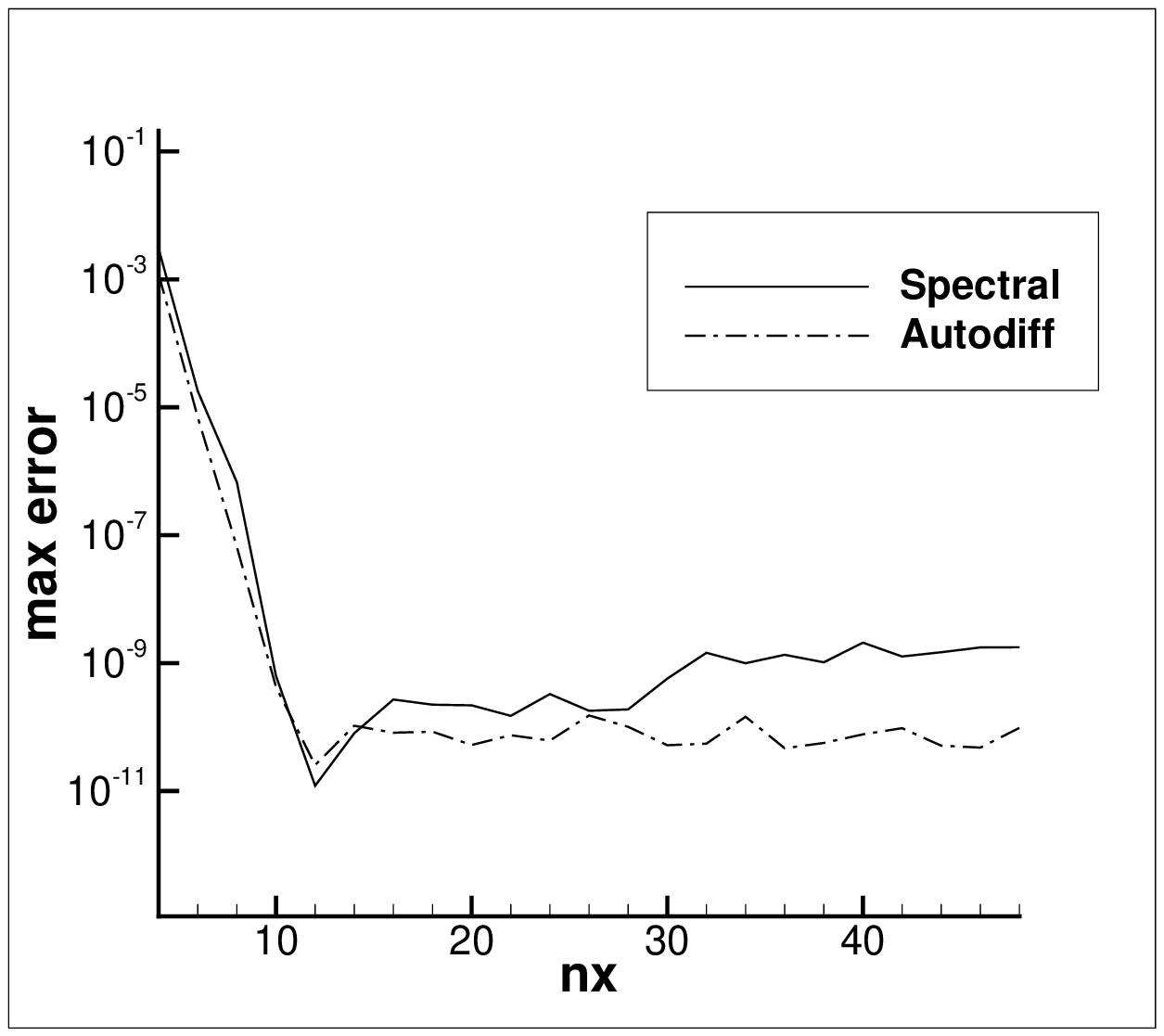}(a)
\includegraphics[width=0.45\textwidth]{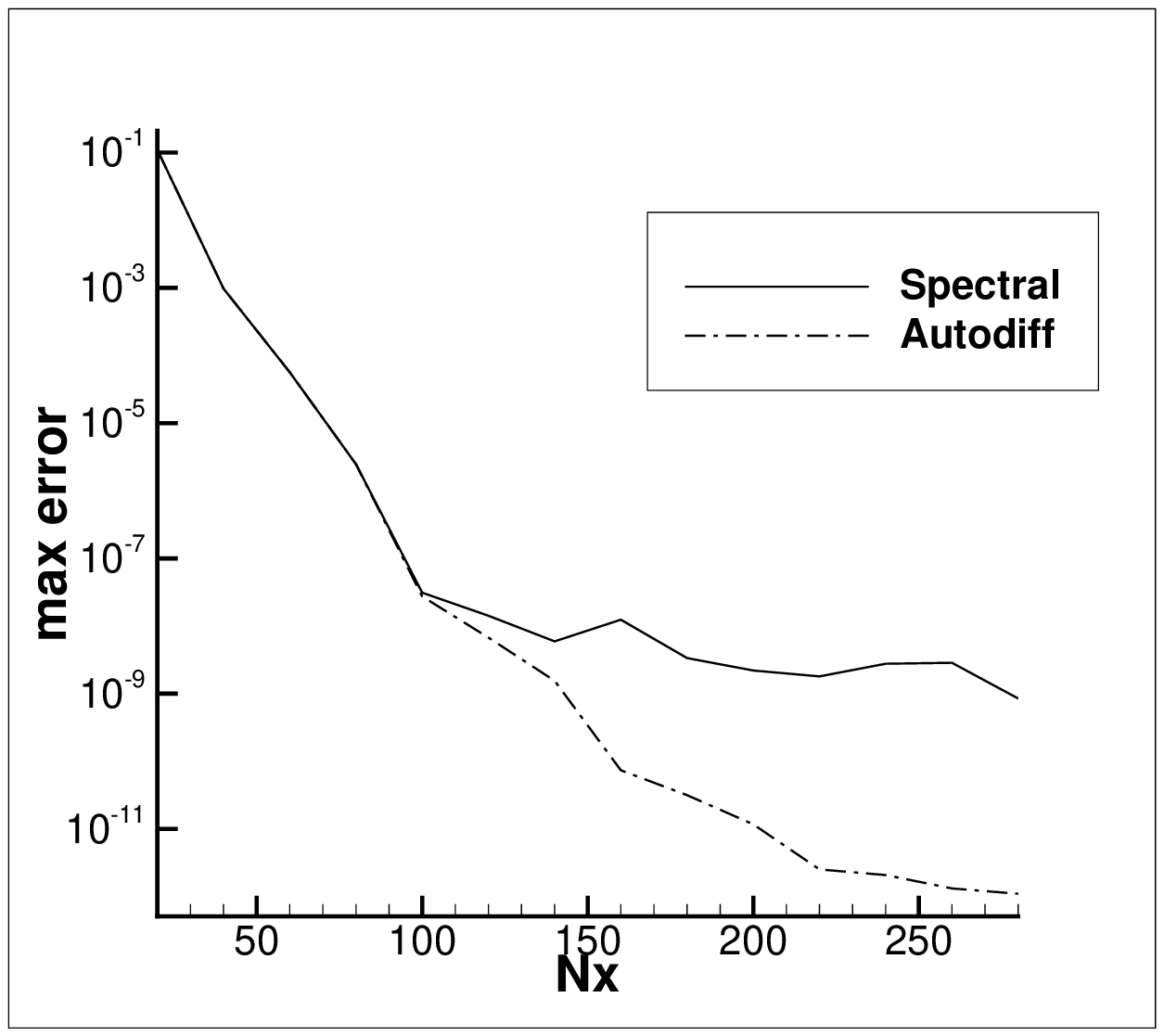}(b)
}
		\caption{Maximum-norm approximation errors from spectral collocation and automatic differentiation versus the number of collocation points $N_x$ (left) and hidden-layer width $M$ (right). Left: $R_m=0.5$ and $M=160$; right: $R_m=1$ and $N_x=40$.
 }\label{fig:nonlinear_poi_compare_max}
	\end{figure}

\begin{figure}[!tb]
	\centerline{
		\includegraphics[width=0.45\textwidth]{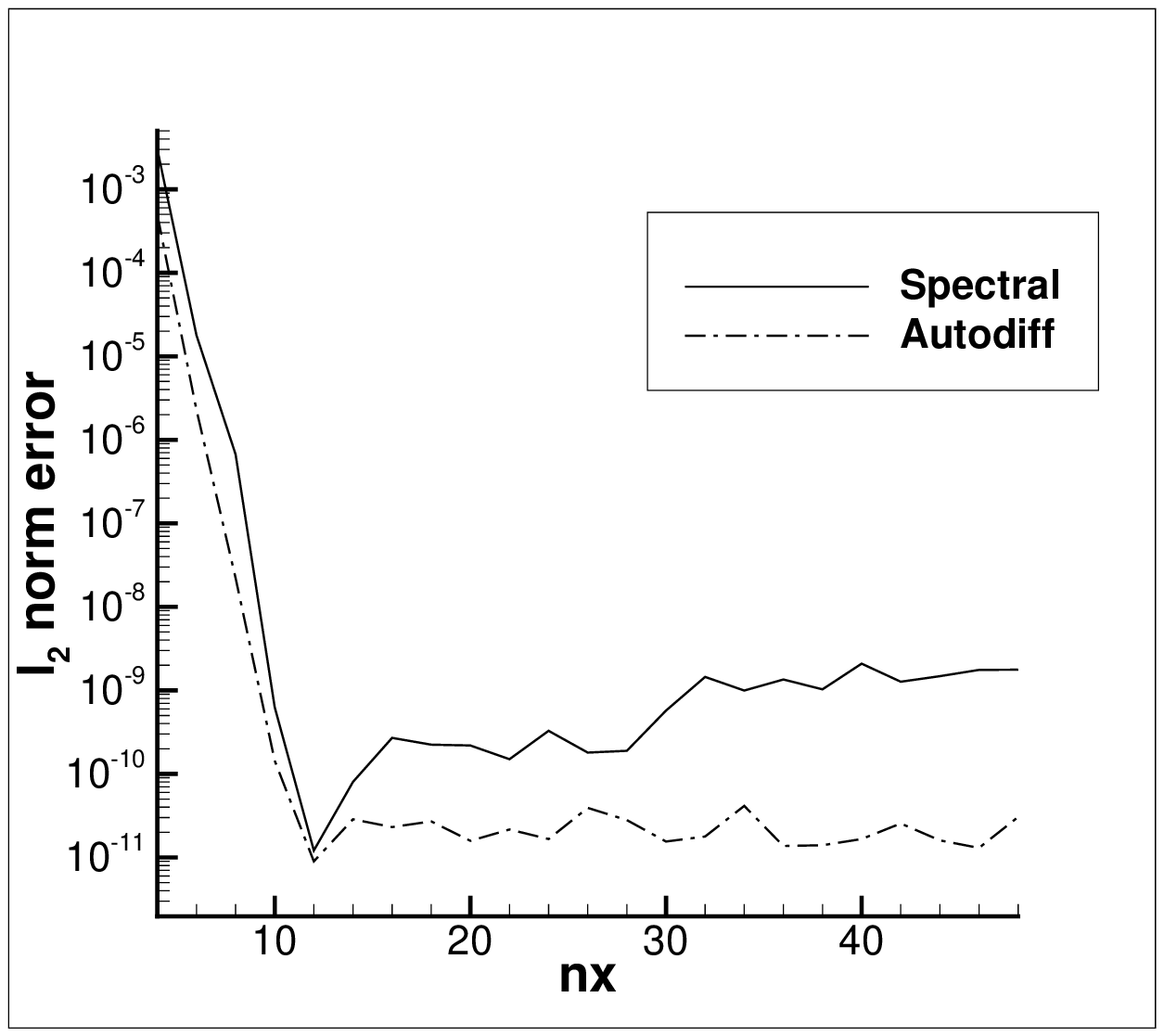}(a)
		\includegraphics[width=0.45\textwidth]{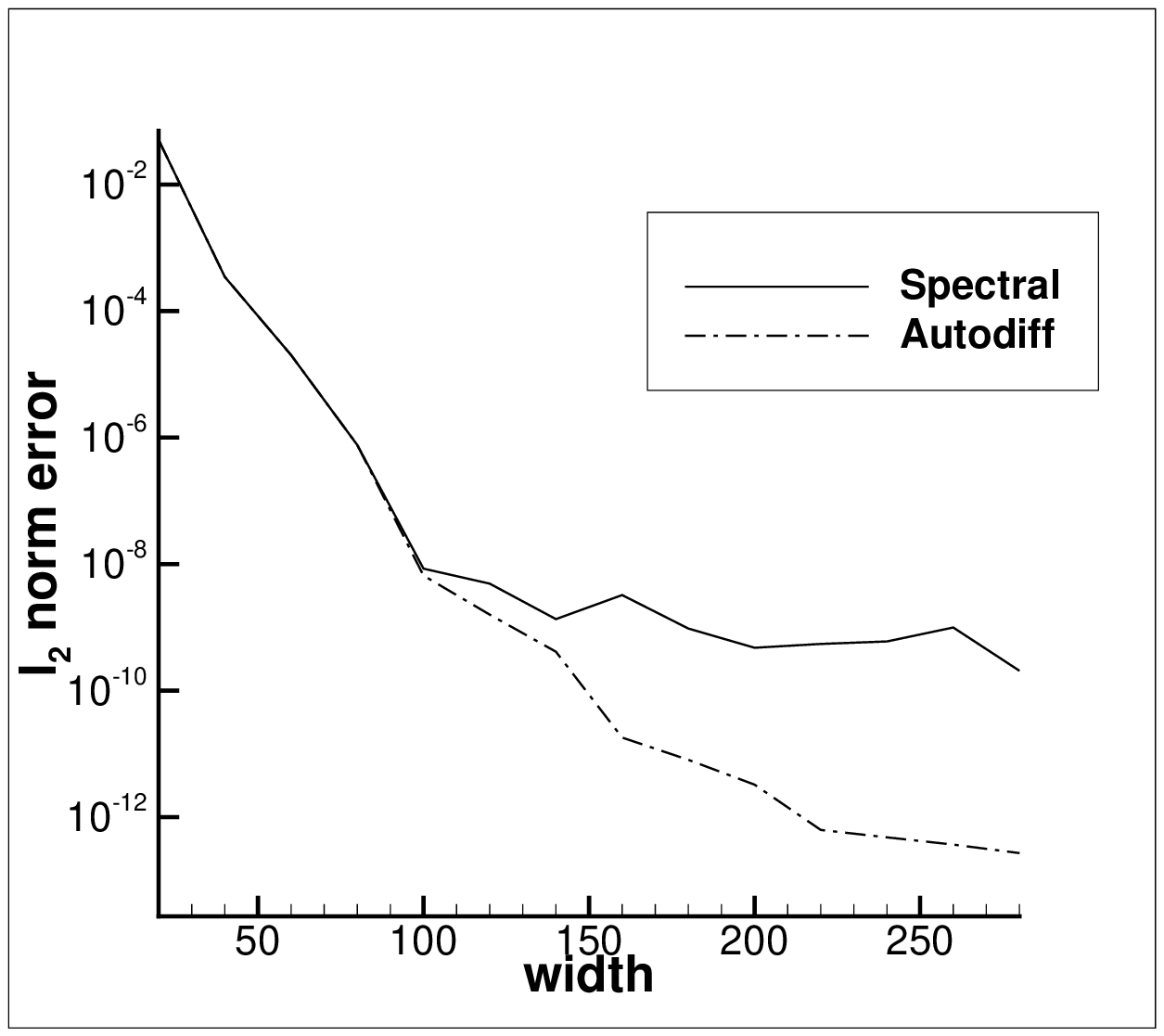}(b)
	}
		\caption{$L^2$-norm approximation errors from spectral collocation and automatic differentiation versus the number of collocation points $N_x$ (left) and hidden-layer width $M$ (right). Left: $R_m=0.5$ and $M=160$; right: $R_m=1$ and $N_x=40$.
		}\label{fig:nonlinear_poi_compare_l2}
\end{figure}

\textbf{Numerical results.} The influence of the random-initialization magnitude $R_m$ on accuracy is presented in Table~\ref{tab:rm_nllinear_poi}. For very small magnitudes ($R_m\leq 0.2$), the error is relatively large ($\sim 10^{-1}$), indicating that the random features are insufficiently expressive to capture the two-dimensional nonlinear profile. Accuracy improves substantially at $R_m=0.4$, with the AD variant reaching an error of $2.70\times 10^{-11}$.
\begin{table}[tb]
	\centering
\begin{tabular}{cccccccc}
	\hline
	$R_m$ & 0.1&0.2&0.4&0.8&1.0&1.2&1.8 \\ \hline 
	$e_{\text{spec}}$ & 4.11E-1& 1.90E-1 & 7.63E-9&1.48E-8&5.52E-8&1.87E-7&2.62E-5\\ 
	$e_{\text{auto}}$ & 4.11E-1&1.90E-1  & 2.70E-11&1.23E-8&5.50E-8&1.96E-7&2.62E-5\\ 
	\hline
\end{tabular}
\caption{Approximation errors for the 2D nonlinear Poisson equation with varying $R_m$ using $M=160$ and $N_x=40$. $e_{\text{spec}}$ and $e_{\text{auto}}$ represent spectral and automatic differentiation maximum errors, respectively.}
\label{tab:rm_nllinear_poi}
\end{table}
The convergence behavior is further explored in Figures~\ref{fig:nonlinear_poi_compare_max} and~\ref{fig:nonlinear_poi_compare_l2}. As shown in panel (a) of Figure~\ref{fig:nonlinear_poi_compare_max}, both differentiation strategies exhibit a steep decay in the maximum error as the number of collocation points $N_x$ increases, eventually plateauing near $N_x\approx 40$. Similarly, panel (b) indicates that increasing the hidden-layer width $M$ is critical for two-dimensional problems; the error decreases consistently as $M$ approaches 300, confirming that a larger feature space is necessary to approximate the two-dimensional nonlinear solution. Consistent with the linear case, the AD variant generally maintains a slight accuracy advantage over spectral collocation across the parameter ranges.

\subsection{Poisson Equation on a Curved Multiply Connected Domain}
\label{sec:wavy_annulus}

We consider the Poisson equation
\begin{equation}
	-\Delta u=f\quad\text{in }\Omega,
	\qquad
	u=g\quad\text{on }\partial\Omega,
\end{equation}
on the curved, multiply connected domain
\begin{equation}
	\Omega=\left\{(r\cos\theta,r\sin\theta):
	r_{\mathrm{in}}(\theta)<r<r_{\mathrm{out}}(\theta),\quad
	0\leq\theta<2\pi\right\},
\end{equation}
where
\begin{equation}
	r_{\mathrm{in}}(\theta)=0.33[1+0.08\cos(3\theta-0.2)],
	\qquad
	r_{\mathrm{out}}(\theta)=1+0.18\cos(5\theta+0.15).
\end{equation}
The domain contains a curved hole and is divided into eight curvilinear
quadrilateral patches. We use the manufactured solution
\begin{equation}
	u(x,y)=\sin(\pi x)\cos(\pi y/2)+0.1\exp(x+y/2),
\end{equation}
and compute the source term and Dirichlet data from this solution.

Each patch contains $13\times13$ CGL points. The local ELM uses tanh features
with fixed parameters sampled uniformly from $[-2,2]$. The feature matrices
are orthonormalized by the local thin QR factorization, and the least-squares
cutoff is $10^{-12}$. We compare the mapped multi-patch Spectral-ELM with the
direct mapped spectral collocation method using the same patch grid and
mapped differential operator. ELM results are reported over three independent
random seeds.

\begin{table}[!htbp]
	\centering
	\small
	\begin{tabular}{lrrrr}
		\hline
		Method & Unknowns & $M_e$ & $\kappa_2$ & $\|u-u_h\|_\infty$ \\
		\hline
		Direct mapped spectral & 1352 & --  & $6.983\mathrm{E}{+2}$ & $1.051\mathrm{E}{-5}$ \\
		Multi-patch Spectral-ELM &  960 & 120 & $2.711\mathrm{E}{+2}$ & $3.933\mathrm{E}{-4}$ \\
		Multi-patch Spectral-ELM & 1120 & 140 & $2.874\mathrm{E}{+2}$ & $7.411\mathrm{E}{-5}$ \\
		Multi-patch Spectral-ELM & 1240 & 155 & $2.881\mathrm{E}{+2}$ & $1.425\mathrm{E}{-5}$ \\
		Multi-patch Spectral-ELM & 1320 & 165 & $2.940\mathrm{E}{+2}$ & $1.082\mathrm{E}{-5}$ \\
		\hline
	\end{tabular}
	\caption{Direct mapped spectral collocation and mapped multi-patch
		Spectral-ELM on the curved domain. ELM errors and condition numbers are medians
		over three random seeds.}
	\label{tab:curved_domain_results}
\end{table}

The direct mapped spectral solution shows rapid convergence. Its maximum error
decreases from $6.63\times10^{-3}$ with $n=7$ to $1.45\times10^{-7}$ with
$n=17$. This direct method is the more accurate choice when the local ELM
width is small. As $M_e$ increases, the Spectral-ELM error decreases steadily.
At $M_e=155$, the median maximum error is $1.43\times10^{-5}$ with 1240
unknowns, compared with $1.05\times10^{-5}$ and 1352 unknowns for direct mapped
spectral collocation. At $M_e=165$, the Spectral-ELM reaches approximately the
same discretization error as the direct method.

The QR basis change is important for numerical stability. For the $n=11$,
$M_e=110$ case, it reduces the scaled condition number from
$1.66\times10^9$ to $1.92\times10^2$ without changing the solution. The result
shows that the mapped multi-patch construction extends CGL differentiation to
a curved domain with a hole. It also shows that the ELM trial space does not
uniformly outperform the direct spectral trial space for a smooth forward
problem.

\begin{figure}[H]
	\centering
	\includegraphics[width=0.98\textwidth]{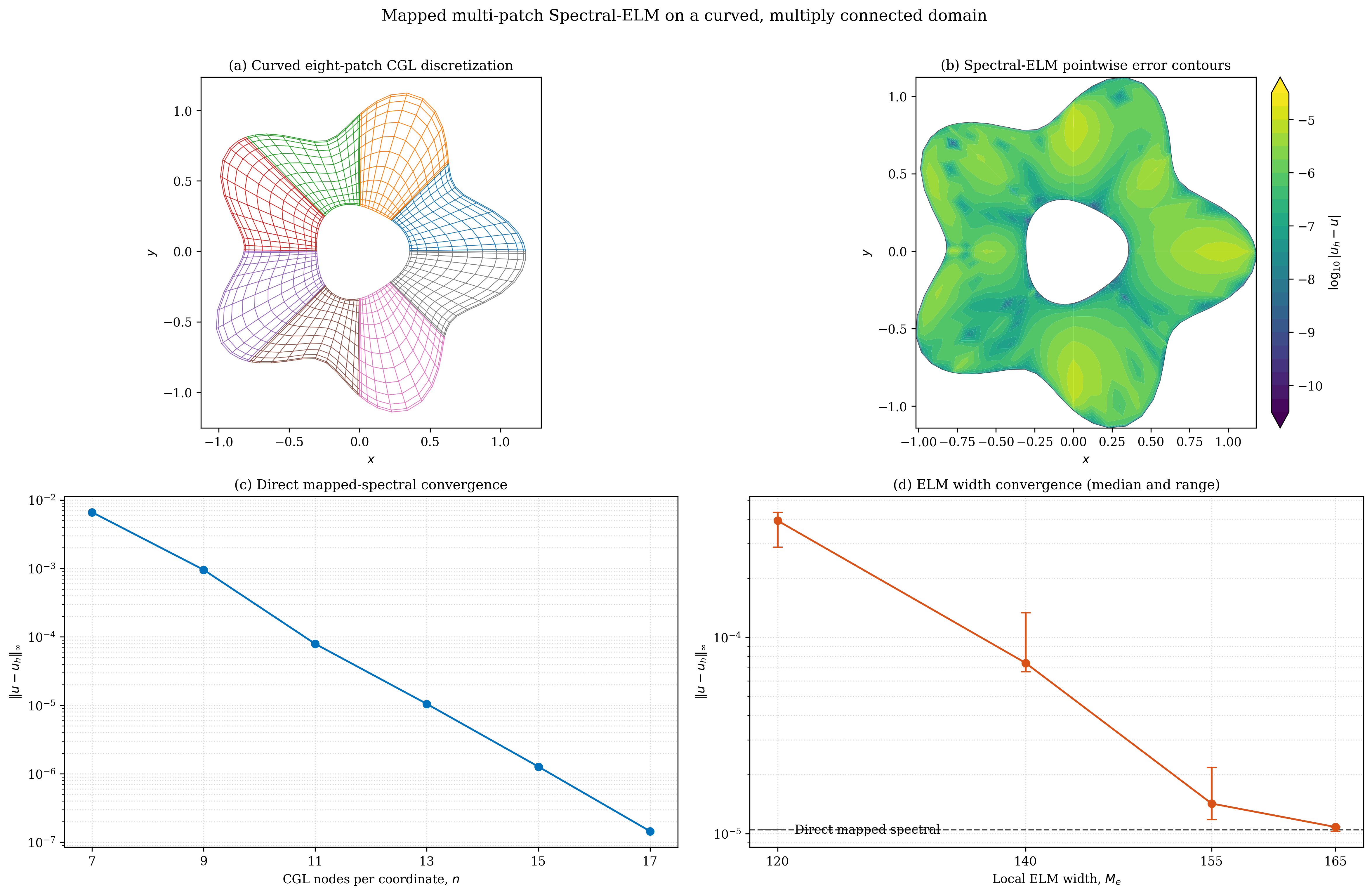}
	\caption{Mapped multi-patch benchmark. (a) Curved patch decomposition and
		local CGL grids. (b) Patchwise filled contours of the pointwise absolute error
		for the best Spectral-ELM run, shown as $\log_{10}|u_h-u|$ with a common color
		scale across all eight patches.
		(c) Convergence of direct mapped spectral collocation. (d) ELM width
		convergence; markers show the median and bars show the range over three seeds.}
	\label{fig:complex_geometry_benchmark}
\end{figure}

\subsection{Burgers' equation}
\textbf{Problem setup.} We further evaluate the ELM framework on the viscous Burgers' equation, a fundamental nonlinear PDE that combines advection and diffusion. We consider the one-dimensional domain $\Omega=[-1,1]$ over the time interval $t\in[0,0.1]$:
\begin{equation*}
    \begin{aligned}
    u_t+uu_x-\mu u_{xx}&=f,\quad x\in \Omega,
    u(-1,t)=u(1,t)&=0,
    u(x,0)&=-\sin(\pi x),
    \end{aligned}
\end{equation*}
where $\mu$ is the viscosity coefficient. We test three regimes: $\mu\in\{1/\pi,1/(10\pi),1/(100\pi)\}$, corresponding to approximate values of $0.318$, $0.0318$, and $0.00318$. Decreasing $\mu$ increases the Reynolds number, leading to sharper gradients that are more challenging to resolve. The source term $f$ is typically zero, and the solution is compared with a reference computed by a high-fidelity viscous Burgers solver. The problem is particularly challenging because steep gradients may form in the solution profile.

\textbf{Special implementation.} The problem is treated as a two-dimensional space--time domain with $m_{\text{in}}=2$. We use a hidden layer of width $M=200$. The collocation grid is constructed using a tensor product of GLJ nodes with $N_x=32$ points in the spatial dimension and $N_t=8$ points in the temporal dimension.

To handle the differential operators, the spectral approach utilizes a Kronecker product assembly to construct the global derivative matrices:
\begin{equation*}
\mathbf{D}_x=\mathbf{I}_t \otimes \mathbf{D}_{x, \text { ori }}, \quad \mathbf{D}_t=\mathbf{D}_{t, \text { ori }} \otimes \mathbf{I}_x, \quad \mathbf{D}_{x x}=\mathbf{I}_t \otimes \mathbf{D}_{x x, \text { ori }}.
\end{equation*}
The boundary conditions and initial conditions are enforced by mapping the respective indices in the space-time grid to a constraint Jacobian matrix. The nonlinear system is solved using the TRF algorithm with a tight tolerance of $1.0\times 10^{-12}$.

\textbf{Numerical results.} The performance of the ELM-AD and ELM-Spectral methods is summarized in Table~\ref{tab:burgers_summary}.

\begin{table}[t]
\centering
\begin{tabular}{cccccc}
\hline
$\mu$ & Variant & $\max|\mathrm{err}|$ & Cost & nfev & Wall [s] \\
\hline
3.18e-01 & AD       & 7.574e-05 & 0.000e+00 & 0 & 0.896 \\
3.18e-01 & Spectral & 1.017e-04 & 0.000e+00 & 0 & 0.441 \\
3.18e-02 & AD       & 6.017e-03 & 0.000e+00 & 0 & 0.631 \\
3.18e-02 & Spectral & 4.841e-03 & 0.000e+00 & 0 & 0.492 \\
3.18e-03 & AD       & 7.601e-03 & 0.000e+00 & 0 & 0.634 \\
3.18e-03 & Spectral & 6.274e-03 & 0.000e+00 & 0 & 0.454 \\
\hline
\end{tabular}
\caption{Comparison of ELM-AD and ELM-Spectral for Burgers' equation on $[-1,1]\times[0,0.1]$ with varying viscosity $\mu$.}
\label{tab:burgers_summary}
\end{table}
The results indicate that both methods provide reasonable approximations for the diffusive case ($\mu \approx 0.3$), with errors between $10^{-5}$ and $10^{-4}$. As the viscosity decreases to $\mu \approx 0.003$, the solution develops steeper gradients, resulting in a degradation of accuracy to approximately $6 \times 10^{-3}$ for the spectral variant and $7 \times 10^{-3}$ for the AD variant. Notably, the spectral collocation method maintains a slight accuracy advantage in the convection-dominated regimes ($\mu \le 0.03$) and consistently requires less computational time, completing the solution in under 0.5 seconds.
\subsection{Steady Allen--Cahn Equation (1D)}
\textbf{Problem setup.} We consider the steady Allen--Cahn equation on the interval $[-L,L]$:
\als
\epsilon^2 u_{xx}+u-u^3=0,
\eals
with Dirichlet boundary data chosen from the exact heteroclinic profile:
\als
u(-L) = \tanh\!\left(\frac{-L}{\sqrt{2}\epsilon}\right), \quad
u(L) = \tanh\!\left(\frac{L}{\sqrt{2}\epsilon}\right).
\eals
We report results for $L=1$ and $\epsilon\in\{1.0,0.5,0.25,0.15,0.10,0.07,0.05\}$. The reference solution on the finite interval is the analytical profile $u_{\text{ref}}(x)=\tanh\!\left(x/(\sqrt{2}\epsilon)\right)$.

\textbf{Neural representation (ELM).} We approximate $u(x)$ using a single-hidden-layer network. The input is the scaled coordinate $\xi=x/L\in[-1,1]$. The hidden layer consists of $M$ $\tanh$ units with preactivation
\als
z_j(\xi;\epsilon)=\frac{s_{\text{base},j}}{\epsilon}(\xi-\mu_j),
\eals
where the centers $\mu_j$ are sampled in $[-1,1]$ and the base slopes $s_{\text{base},j}$ are sampled log-uniformly in $[0.1,10]$, meaning that $\log_{10}(s_{\text{base},j})$ is uniformly sampled from $[-1,1]$. For each $\epsilon$, only the output-layer coefficients $\beta\in\bbR^M$ are optimized; the hidden parameters are randomly assigned and fixed. This $\epsilon$-dependent scaling of the hidden slopes makes the dictionary compatible with the $\epsilon$-scale interface layer in the Allen--Cahn profile.

\textbf{Collocation grid and two differentiation strategies.} We place $n_x$ Chebyshev--Gauss--Lobatto nodes on $x\in[-L,L]$ and consider two ways to form $u_{xx}$ from the hidden-layer features $v(x)\in\bbR^{n_x\times M}$. The first strategy uses automatic differentiation (AD). We compute $d^2v/d\xi^2$ by forward-mode AD in $\xi$ and map it to $x$ through $u_{xx}=(1/L^2)u_{\xi\xi}$. The second strategy uses spectral collocation: we compute $v_{xx}$ with the Chebyshev $D_{xx}$ matrix on the same collocation grid. We solve a nonlinear least-squares problem for the output coefficients $\beta$ by minimizing $\|r\|_2$, where
\als
r = [w^{1/2} (\epsilon^2 (v_{xx} \beta)+v\beta)-(v\beta)^3; w_{\bc}(v_{\bc}\beta-u_{\bc})],
\eals
with $w$ denoting the Clenshaw--Curtis quadrature weights on the PDE rows and $w_{\bc}$ the Dirichlet weight. To improve conditioning, all feature blocks are column-scaled so that each column has unit RMS over the stacked PDE and boundary-condition rows; the Jacobian is right-scaled consistently. We perform $\epsilon$ continuation from a larger, easier value to a smaller, stiffer value, reusing the previous solution as an initializer by projecting it onto the current feature space through a least-squares warm start. This procedure greatly reduces the number of nonlinear iterations and avoids poor local minima. The output coefficients are found with a nonlinear least-squares solver from SciPy. We provide the exact Jacobian of the stacked residual and use tight stopping tolerances (about $10^{-12}$), with a maximum of 400 function evaluations. The solver is run separately for the AD and spectral variants on the same CGL grid. Continuation proceeds from easy to stiff regimes according to the schedule $\epsilon\in\{1.0,0.5,0.25,0.15,0.10,0.07,0.05\}$.

At each $\epsilon$, we report the maximum-norm error $e=\max_k|u_{\text{approx}}(x_k)-u_{\text{ref}}(x_k)|$, final least-squares cost, number of function evaluations, and wall-clock time. The AD variant maintains near-machine precision for $\epsilon\geq 0.25$ and degrades smoothly as the interface sharpens; the ELM spectral variant typically achieves errors from $10^{-11}$ to $10^{-5}$ across the tested range, with competitive or lower cost. A representative run produced the summary in Table~\ref{tab:ac_elm_summary}.

\begin{table}[t]
\centering
\begin{tabular}{cccccc}
\hline
$\varepsilon$ & Variant   & $\max|\mathrm{err}|$ & Cost        & nfev & Wall [s] \\
\hline
1.00  & AD        & 1.221e-15 & 1.897e-29 & 19 & 0.798 \\
1.00  & Spectral  & 2.186e-11 & 5.870e-26 & 34 & 1.249 \\
0.50  & AD        & 1.066e-14 & 2.816e-28 & 28 & 1.361 \\
0.50  & Spectral  & 4.836e-12 & 3.963e-28 & 20 & 0.853 \\
0.25  & AD        & 1.967e-13 & 2.173e-28 & 33 & 1.444 \\
0.25  & Spectral  & 3.141e-10 & 1.403e-28 & 23 & 1.019 \\
0.15  & AD        & 2.106e-03 & 1.731e-28 & 23 & 0.934 \\
0.15  & Spectral  & 1.428e-07 & 2.057e-23 &  6 & 0.434 \\
0.10  & AD        & 7.414e-03 & 7.651e-31 & 26 & 0.963 \\
0.10  & Spectral  & 1.233e-05 & 2.884e-24 & 20 & 0.541 \\
0.07  & AD        & 1.586e-02 & 8.192e-28 & 31 & 1.240 \\
0.07  & Spectral  & 1.793e-05 & 3.848e-25 & 16 & 0.363 \\
0.05  & AD        & 2.057e-02 & 1.694e-25 & 29 & 1.314 \\
0.05  & Spectral  & 2.564e-05 & 2.275e-26 & 17 & 0.320 \\
\hline
\end{tabular}
\caption{Steady Allen--Cahn boundary-value problem on $[-L,L]$ with Dirichlet data, comparing ELM+AD and ELM+Spectral (without an interior anchor). Continuation uses $\varepsilon\in\{1.0,0.5,0.25,0.15,0.10,0.07,0.05\}$ on the same CGL grid. Nonlinear least squares uses TRF with $g_{\mathrm{tol}}=f_{\mathrm{tol}}=x_{\mathrm{tol}}=10^{-12}$ and a maximum of 400 function evaluations. Reported quantities are the maximum error $e=\max|u_{\text{approx}}-u_{\text{ref}}|$, final cost, number of function evaluations (nfev), and wall time (s).}
\label{tab:ac_elm_summary}
\end{table}
These trends are consistent with increasing stiffness as $\epsilon$ decreases and illustrate that (i) both discretizations benefit from feature scaling and continuation and (ii) spectral differentiation paired with quadrature-weighted residuals and column scaling remains robust for small $\epsilon$.
\begin{figure}[t]
  \centering
  \begin{subfigure}[t]{0.32\linewidth}
  \centering
    \includegraphics[width=\linewidth]{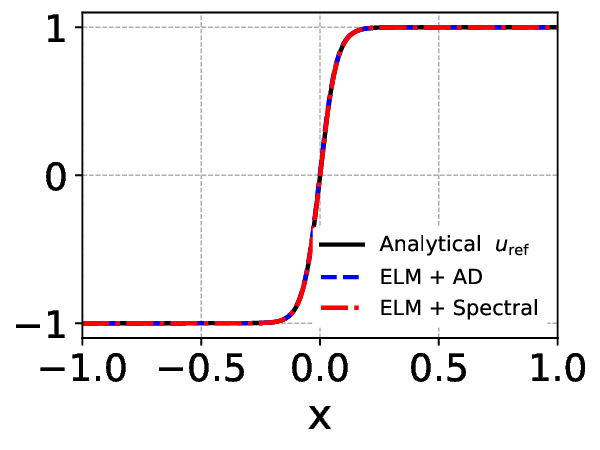} 
    \caption{Analytical reference and ELM solutions (AD and spectral).}
  \end{subfigure}\hfill
  \begin{subfigure}[t]{0.32\linewidth}
    \centering
    \includegraphics[width=\linewidth]{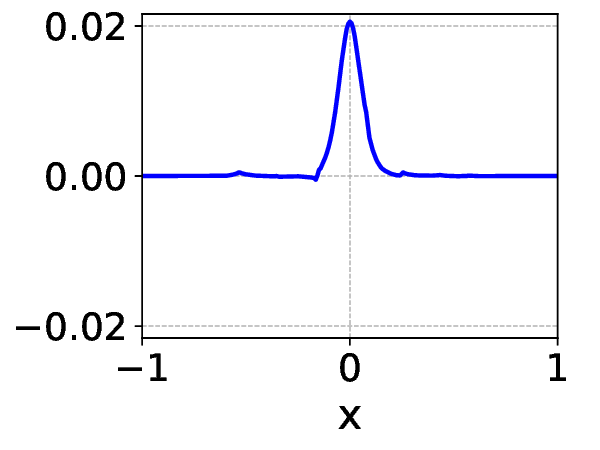}
    \caption{Error of ELM+AD.}
  \end{subfigure}\hfill
  \begin{subfigure}[t]{0.32\linewidth}
    \centering
    \includegraphics[width=\linewidth]{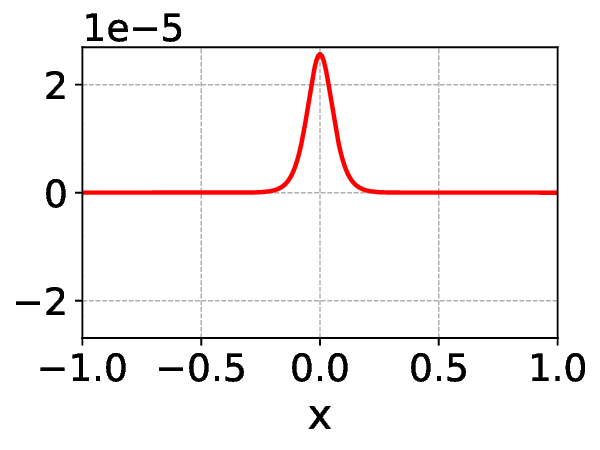}
    \caption{Error of ELM+Spectral.}
  \end{subfigure}

    \caption{Comparison of ELM+AD and ELM+Spectral solutions for the steady Allen--Cahn equation at $\varepsilon=0.05$.
    (a) Reference and neural solutions, (b)--(c) pointwise error distributions.}
    \label{fig:ac_elm_compare}
\end{figure}

The Allen--Cahn test shown in Figure~\ref{fig:ac_elm_compare} visualizes the accuracy of the ELM-based solvers at the smallest continuation step, $\epsilon=0.05$. Panel (a) compares the analytical reference profile $u_{\text{ref}}$ with the two neural solutions obtained using ELM+AD and ELM+Spectral collocation. Both neural approximations reproduce the sharp interface structure and saturation plateaus with relatively good agreement. Panels (b) and (c) display the corresponding pointwise errors. The ELM+AD solution shows increasing deviation near the interface as $\epsilon$ decreases, while the ELM+Spectral solution maintains uniformly small errors across the domain. The maximum absolute errors, $2.1\times 10^{-2}$ and $2.6\times 10^{-5}$, respectively, are consistent with the quantitative results in Table~\ref{tab:ac_elm_summary}.
\subsection{KdV equation: matched derivative comparison}
\label{sec:matched_kdv}

We solve the Korteweg--de Vries equation
\begin{equation}
	u_t+6uu_x+u_{xxx}=0,\qquad (x,t)\in[-20,20]\times[0,1.5],
\end{equation}
with the single-soliton solution
\begin{equation}
	u_{\rm ref}(x,t)=\frac{c}{2}\operatorname{sech}^2
	\left[\frac{\sqrt c}{2}(x-ct-x_0)\right],
	\qquad c=1,\quad x_0=-5.
\end{equation}

The earlier comparison used different trial functions and different collocation points. It could not isolate the effect of differentiation. The revised comparison uses the same trial function in both cases,
\begin{equation}
	u_M(x,t)=u_{\rm ref}(x,0)+t\left(1-\frac{x^2}{L^2}\right)
	\sum_{m=1}^{M}\beta_m\tanh(w_{xm}x+w_{tm}t+b_m).
\end{equation}
Thus the initial condition is exact and the small boundary tails are fixed in the same way. Both variants use the same $129\times41$ CGL points, $M=220$ features, column scaling, zero initial coefficient vector, residual, Jacobian formula, trust-region least-squares solver, tolerances, and off-grid evaluation mesh. The only change is this: one variant applies CGL matrices to the nodal trial function, while the other uses analytic derivatives of the same tanh features. We use three fixed seeds.

\begin{table}[H]
	\centering
	\small
	\begin{tabular}{lrrrr}
		\hline
		Derivative calculation & Median max. error & Median RMS error & Median rel. $L^2(T)$ & Median nfev \\
		\hline
		CGL matrices & $1.97\times10^{-2}$ & $2.93\times10^{-3}$ & $2.67\times10^{-2}$ & 19 \\
		Analytic formulas & $1.78\times10^{-2}$ & $2.66\times10^{-3}$ & $2.39\times10^{-2}$ & 19 \\
		\hline
	\end{tabular}
	\caption{Matched KdV comparison. Every part of the discretization and solve is held fixed except the derivative calculation. Medians are over three seeds.}
	\label{tab:kdv_matched}
\end{table}

\begin{figure}[H]
	\centering
	\begin{subfigure}[t]{0.45\linewidth}
		\centering
		\includegraphics[width=\linewidth]{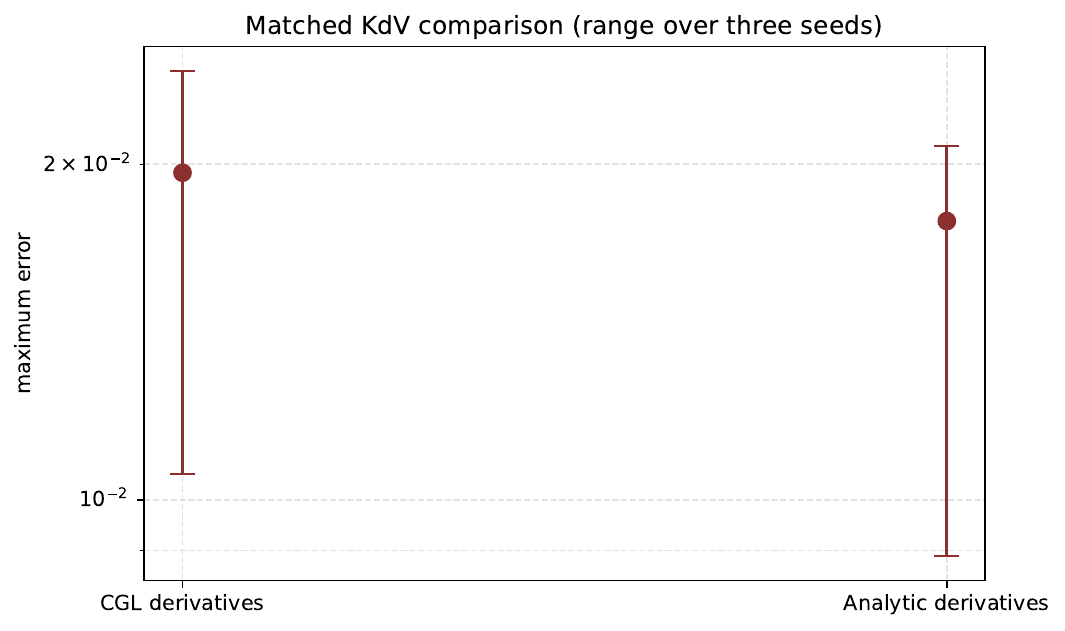}
		\caption{Maximum error over three seeds.}
	\end{subfigure}\hfill
	\begin{subfigure}[t]{0.52\linewidth}
		\centering
		\includegraphics[width=\linewidth]{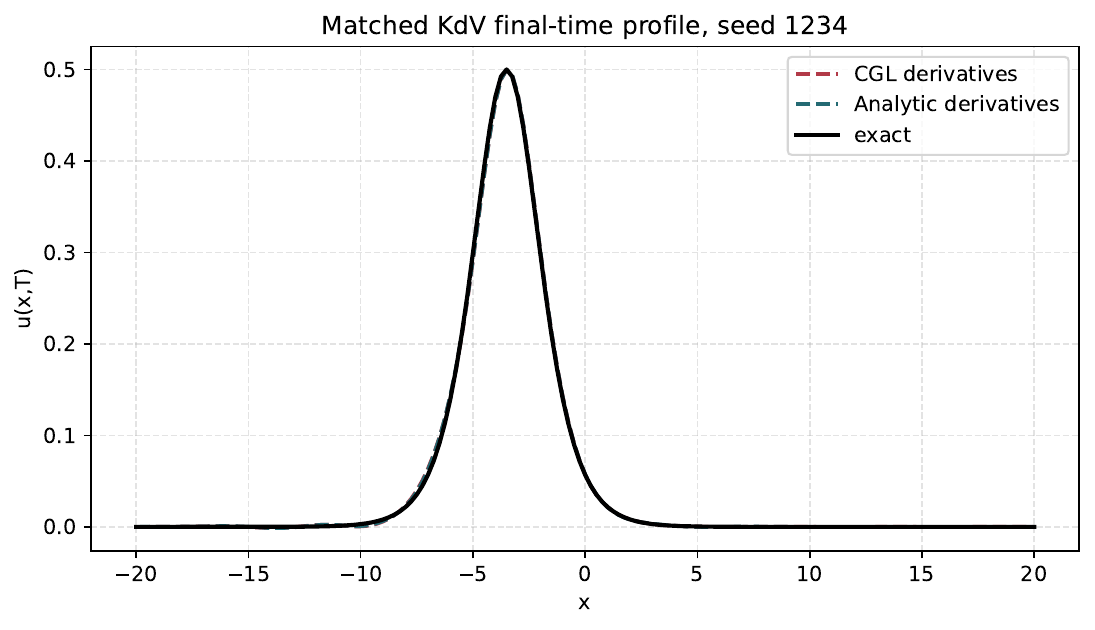}
		\caption{Final-time profile for seed 1234.}
	\end{subfigure}
	\caption{Matched KdV derivative comparison.}
	\label{fig:kdv_matched}
\end{figure}

The two errors are close. Analytic differentiation is slightly more accurate at this grid resolution, but the ranges overlap. Therefore, the earlier large difference was mainly caused by changing the trial function and residual sampling, not only by changing the derivative formula. The result also gives a practical requirement for third derivatives: the CGL grid must resolve every random feature. If the feature weights create spatial scales that are too small for the grid, applying $D_x^3$ amplifies the unresolved part. We use a spatial feature scale of 0.30 here; a coarser $65$-point grid did not resolve the third derivatives reliably and was rejected before making the comparison.

\subsection*{Overall numerical findings}

The controlled two-dimensional elliptic test gives the clearest comparison. Direct CGL reaches a maximum error of $3.09\times10^{-10}$, while Spectral-ELM reaches a median of $4.03\times10^{-8}$ with fewer unknowns. The direct method is therefore preferable for this smooth square-domain forward problem when maximum accuracy is the main goal. The fixed-feature formulation is not a replacement for the complete nodal polynomial space.

The same control separates feature choice from derivative choice. Spectral-ELM and analytic ELM use identical features and give errors of the same order. The global RFM baseline has a much larger condition number and error in this setting. The domain-scaled TransNet sampling rule improves the feature space, but the resulting error is still above the two ELM controls. These rows are limited controls, not complete reproductions of localized RFM or auxiliary-trained TransNet. SNN has a separate feature-training stage and is therefore compared at the method level rather than mixed into the fixed-feature timing and error table.

The matched KdV test leads to the same main conclusion. When the trial function, nodes, features, and nonlinear solver are fixed, CGL and analytic derivatives have close errors. Analytic differentiation is slightly better at the tested resolution. The third derivative also makes the resolution condition explicit: the CGL grid must resolve the smallest scale present in every random feature.

The curved-domain Poisson test addresses a different question. Eight curved quadrilateral patches represent a wavy annulus without embedding it in a box. At $n=13$ points per reference direction, direct mapped spectral collocation gives a maximum error of $1.05\times10^{-5}$ with 1352 unknowns. The local ELM with 165 features per patch gives a median maximum error of $1.08\times10^{-5}$ with 1320 unknowns. The direct method converges from $6.63\times10^{-3}$ at $n=7$ to $1.45\times10^{-7}$ at $n=17$. Thus the mapping and interface equations handle the geometry, while the choice between a nodal and a fixed-feature trial space remains problem dependent.

The one-dimensional Poisson, nonlinear Poisson, Burgers, and Allen--Cahn tests show how feature scale, stiffness, and random initialization affect the two ELM variants. They should be read as equation-specific tests, not as evidence that one derivative mechanism always wins. Taken together, the results support a limited statement: CGL matrices provide a practical derivative operator for nodal values of a fixed-feature approximation, provided that the grid resolves those values; analytic derivatives remain equally valid and can be better for high-order operators.

\section{Concluding Remarks}

We have developed a mapped multi-patch Spectral-ELM framework for the numerical solution of smooth PDEs in low-dimensional domains. The method uses a fixed-feature ELM representation, with Chebyshev--Gauss--Lobatto differentiation applied on reference grids. Curved physical domains are decomposed into quadrilateral patches, each obtained by mapping a reference square into the physical domain. The corresponding metric terms transform the reference derivatives into physical derivatives, while neighboring patches are coupled by enforcing continuity of the solution and its physical normal flux. A local QR construction provides a more stable basis for the sampled feature space and reduces the effects of near-linear dependence among the discrete features.

The numerical study clarifies the role of the proposed method relative to both classical spectral collocation and analytically differentiated fixed-feature models. Direct CGL collocation remains the most accurate and straightforward choice for smooth forward problems on simple geometries when the full nodal representation is computationally practical. Spectral-ELM instead restricts the nodal approximation to a lower-dimensional feature space, which can reduce the number of unknown coefficients but introduces dependence on the quality and random realization of the features. The results therefore show that the feature space, rather than the differentiation rule alone, is often the main factor controlling the final approximation error.

When the trial function and numerical solver are held fixed, spectral and analytic differentiation produce comparable errors in the examples considered. In the matched KdV experiment, analytic differentiation is slightly more accurate at the tested resolution, while CGL differentiation offers a unified nodal procedure for computing derivatives of different orders. On the curved-domain problem, mapped Spectral-ELM and mapped direct CGL collocation achieve similar accuracy at comparable numbers of unknowns. These results suggest that the main value of Spectral-ELM lies not in replacing direct spectral collocation on standard domains, but in combining a compact local feature representation with established high-order differentiation and a consistent treatment of mapped multi-patch geometry.

The present formulation is intended for smooth PDEs in one to three dimensions. Although the number of feature coefficients may be smaller than the number of collocation points, the trial function must still be evaluated on a tensor-product CGL grid, whose size grows as $n^d$ when $n$ points are used in each of $d$ coordinate directions. The three-dimensional example confirms that the formulation extends naturally beyond two dimensions, but it is not intended to demonstrate scalability for genuinely high-dimensional problems. Similarly, the current experiments do not address nonsmooth solutions or geometries obtained directly from general computer-aided design descriptions.

Several extensions follow naturally from the present work. Automatic generation and quality control of curved patches would make the geometric construction easier to apply to more complicated domains. Nonconforming interfaces and three-dimensional mapped elements would broaden the range of admissible discretizations. Rank-revealing feature selection could provide a more systematic way to construct compact local trial spaces, while suitable preconditioners would improve the solution of increasingly large patchwise systems. For higher-dimensional PDEs, combining fixed-feature approximation with sparse grids, low-rank tensor formats, or dimension-adaptive nodal sets offers a promising direction for overcoming the growth of the underlying tensor grid.

\bibliographystyle{plain}
\bibliography{references}
\end{document}